
\let\nobibtex = t
\let\noarrow = t

\let\indexproofingtrue = t
\input eplain
\let\indexproofingtrue = t
\input epsf.tex
\input rotate.tex

\def\mystrut{\vrule width 0pt height 6pt depth 3pt}

\def \punclist{\vskip -14pt\vbox{\hbox{\vrule height .8mm width .8mm}}}
\def \puncrem{\vskip 2pt\vbox{\hbox{\vrule height .8mm width .8mm}}\bigskip}

\definecontentsfile{tables}
\definecontentsfile{figs}
\catcode`\;=\active \def ;{\unskip\kern.2em\string;\ }
\catcode`\:=\active \def :{\unskip\kern.25em\string:\ \ignorespaces}
\pretolerance 500 \tolerance 1000 \brokenpenalty 5000

\frenchspacing \hoffset -.4mm \hsize 16.2cm \vsize 23cm
\parindent 0cm


\catcode`\Ž=\active \def Ž{\'e} \catcode`\'=\active \def '{\"e}
\catcode`\=\active \def {\`e} \catcode`\ˆ=\active \def ˆ{\`a}
\catcode`\=\active \def {\^e} \catcode`\=\active \def {\c c}
\catcode`\‰=\active \def ‰{\^a} \catcode`\Š=\active \def Š{\"a}
\catcode`\=\active \def {\`u} \catcode`\ž=\active \def ž{\^u} 
\catcode`\"=\active \def "{\^\i} \catcode`\•=\active \def •{\"\i} 
\catcode`\™=\active \def ™{\^o} \catcode`\š=\active \def š{\"o}
 \catcode`\Ÿ=\active \def Ÿ{\"u}

\font \mathbb=msbm10
\font \mathbbs=msbm8

\def \fix(#1,#2){#1^{\star#2}}

\def\det{\mathop{\hbox{\rm det}}}

\def\sup{\mathop{\hbox{\rm sup}}}
\def\max{\mathop{\hbox{\rm max}}}

\def \Inm{\big\{(i^1,\ldots, i^m)\colon i^\ell\in\N,\, 1\leq i^\ell\leq n,\, 
i^{\ell_1}\neq i^{\ell_2}\, \hbox{ if } \ell_1\neq \ell_2\big\}}

\def \Onm{\big\{(i^1,\ldots, i^m)\colon i^\ell\in\N,\, 1\leq i^\ell\leq n\big\}}
\def \InminOnm{\big\{i\in\Omega_n^m\colon \ell_1\neq \ell_2  
\Rightarrow i^{\ell_1}\neq i^{\ell_2}\big\}}

\def \Pnfa {P^n_{f_A}}

\def \Enfa {E^n_{f_A}}

\def\ind{\hbox{\mathbb I}}

\def \indic#1{\,\ind\left\{#1\right\}}
\def \indicsmall#1{\,\ind_{\left\{#1\right\}}}

\def \smallZ{\hbox{\mathbbs Z}}
\def \smallR{\hbox{\mathbbs R}}

\def \R{\hbox{\mathbb R}}
\def \Z{\hbox{\mathbb Z}}
\def \N{\hbox{\mathbb N}}
\def \smallR{\hbox{\mathbb R}}

\parskip \baselineskip
\parindent 0cm
\topmargin=1.2in
\bottommargin=1in
\leftmargin=1.25in
\rightmargin=1.25in

\newcount\notitre 
\newcount\nosection
\newcount \noproposition
\newcount \nolemma
\newcount \nocorollary
\newcount \nodefinition

\font \bigtitlefont = cmbx12 at 14pt
\font \titlefont = cmbx12
\font \sectionfont = cmbx10 at 11pt

\font \abstracttitle = cmbx10 at 11pt
\font \abstract = cmr10

\font \rm = cmr10 at 11pt

\font \smallsf = cmss8
\font \tt = cmtt10
\font \smalltt = cmtt8

\def \titre  #1\par%
{ \bigskip\penalty -2000\nosection=0 \advance\notitre by 1
{\parindent=0 cm \titlefont \the\notitre. #1}\bigskip
\writenumberedtocentry{chapter}{#1}{\the\notitre}}

\def \freesection  #1\par%
{ \bigskip\penalty -2000
{\parindent=0 cm \titlefont  #1}\bigskip
\writetocentry{freesection}{#1}}

\def \section #1\par%
{\bigskip\advance\nosection by 1\bigskip%
{\parindent=0 cm \sectionfont \the\notitre.\the\nosection\ #1}\par
\writenumberedtocentry{section}{#1}{\the\nosection}
}

\def \subsection #1\par{\bigskip { \bf #1}\par
\writetocentry{subsection}{ #1}}

\footline{\abstract \ifnum\pageno=1\else\hfil\the\pageno\hfil\fi}
\rm
\def \spec#1{{\left\langle\,#1\,\right\rangle}}

\def \proposition #1\par{
\advance \noproposition by 1
\bigskip{\bf Proposition \the\notitre.\the\noproposition \ \ (#1)}\par
\writetocentry{subsection}{ #1}
}

\def \lemma #1\par{
\advance \nolemma by 1
\smallskip{\bf Lemma \the\notitre.\the\nolemma \ \ (#1)}\par
\writetocentry{subsection}{#1}
}
\def \proof{\bigskip}

\def \corollary #1\par{
\advance \nocorollary by 1
\smallskip{\bf Corollary \the\notitre.\the\nocorollary \ \ (#1)}\par
\writetocentry{subsection}{ #1}
}

%

\vglue .5in
{\bigtitlefont Independent Component Analysis and estimation of a quadratic functional}
\bigskip 
\bigskip

\leftline{{\bf Pascal Barbedor}} 
{\it
\leftline{Laboratoire de probabilit\'es et mod\'eles al\'eatoires}
\leftline{CNRS-UMR 7599, Universit\'e Paris VI et 
Universit\'e Paris VII}
\leftline{PO box 7012, 75251 PARIS Cedex 05 (FRANCE)}
\leftline{barbedor@proba.jussieu.fr}

} \bigskip
 
\bigskip
 
 \centerline{\abstracttitle Abstract}
{\abstract \advance \leftskip by .25in
\advance \rightskip by .25in 
Independent component analysis (ICA) is linked up with the problem of estimating a non linear functional of a density for which optimal estimators are well known.
The precision of ICA is analyzed from the viewpoint of functional spaces in the wavelet framework.
In particular, it is shown that, under Besov smoothness conditions,
parametric rate of convergence is  achieved by a U-statistic estimator of the wavelet ICA contrast,
while the previously introduced plug-in estimator $\hat C^2_j$, with moderate computational cost, has a rate in $n^{-4s\over 4s+d}$.

 \par
}
\bigskip


\def \smallZ{\hbox{\mathbbs Z}}
\def \smallR{\hbox{\mathbbs R}}

\def\sup{\mathop{\hbox{\rm sup}}}

\def \det{\mathop{\hbox{\rm det}}}

\def \br{\hfill\break }
\def \Np (#1){\|#1\|_p}
\def \Ndeux (#1){\|#1\|_2}

\def \ie{{\it i.e.\ }}
\def \R{\hbox{\mathbb R}}
\def \Z{\hbox{\mathbb Z}}
\def \N{\hbox{\mathbb N}}
\def \echn(#1){#1_1\ldots,#1_n}
\def \echd(#1){{#1^1\!\!\ldots,#1^d}}

\def\max {\mathop{\hbox{\rm max}}}

\def \prodalpha {\alpha_{jk^1}\ldots\alpha_{jk^d}}

\def \pralpha #1,#2{\alpha_{j#1^1}^#2\ldots\alpha_{j#1^d}^#2}
\def \hatpralpha #1,#2{\hat\alpha_{j#1^1}^#2\ldots\hat\alpha_{j#1^d}^#2}

\def \hatprodalpha {\hat\alpha_{jk^1}\ldots\hat\alpha_{jk^d}}

\def \mulalpha {\alpha_{j\echd(k)}}

\def \Pnfa {P^n_{f_A}}

\def \Enfa {E^n_{f_A}}

\def \Pfa {P_{f_A}}

\def \punclist{\vskip -14pt\vbox{\hbox{\vrule height .8mm width .8mm}}}
\def \puncrem{\vskip 2pt\vbox{\hbox{\vrule height .8mm width .8mm}}\bigskip}
\def \endproof{
\vbox{\hrule\hbox to 5pt{\vrule height 4,2 pt\hfil\vrule}\hrule
}
\bigskip}


\rm


\def \smallZ{\hbox{\mathbbs Z}}
\def \smallR{\hbox{\mathbbs R}}

\def\ind{\hbox{\mathbb I}}

\def\sup{\mathop{\hbox{\rm sup}}}

\def \det{\mathop{\hbox{\rm det}}}

\def \fix(#1,#2){#1^{\star#2}}

\def \br{\hfill\break }
\def \Np (#1){\|#1\|_p}
\def \Ndeux (#1){\|#1\|_2}

\def \ie{{\it i.e.\ }}
\def \R{\hbox{\mathbb R}}
\def \Z{\hbox{\mathbb Z}}
\def \N{\hbox{\mathbb N}}
\def \echn(#1){#1_1\ldots,#1_n}
\def \echd(#1){{#1^1\!\!\ldots,#1^d}}

\def\max {\mathop{\hbox{\rm max}}}

\def\mystrut{\vrule width 0pt height 6pt depth 3pt}

\def \Inm{\big\{(i^1,\ldots, i^m)\colon i^\ell\in\N,\, 1\leq i^\ell\leq n,\, 
i^{\ell_1}\neq i^{\ell_2}\, \hbox{ if } \ell_1\neq \ell_2\big\}}

\def \Onm{\big\{(i^1,\ldots, i^m)\colon i^\ell\in\N,\, 1\leq i^\ell\leq n\big\}}
\def \InminOnm{\big\{i\in\Omega_n^m\colon \ell_1\neq \ell_2  
\Rightarrow i^{\ell_1}\neq i^{\ell_2}\big\}}

\def \indic#1{\,\ind\left\{#1\right\}}
\def \indicsmall#1{\,\ind_{\left\{#1\right\}}}

\def \spec#1{{\left\langle\,#1\,\right\rangle}}

\def \prodalpha {\alpha_{jk^1}\ldots\alpha_{jk^d}}

\def \pralpha #1,#2{\alpha_{j#1^1}^#2\ldots\alpha_{j#1^d}^#2}
\def \hatpralpha #1,#2{\hat\alpha_{j#1^1}^#2\ldots\hat\alpha_{j#1^d}^#2}

\def \hatprodalpha {\hat\alpha_{jk^1}\ldots\hat\alpha_{jk^d}}

\def \mulalpha {\alpha_{j\echd(k)}}

\def \Pnfa {P^n_{f_A}}

\def \Enfa {E^n_{f_A}}

\def \Pfa {P_{f_A}}

\def \punclist{\vskip -14pt\vbox{\hbox{\vrule height .8mm width .8mm}}}
\def \puncrem{\vskip 2pt\vbox{\hbox{\vrule height .8mm width .8mm}}\bigskip}
\def \endproof{
\vbox{\hrule\hbox to 5pt{\vrule height 4,2 pt\hfil\vrule}\hrule
}
\bigskip}

\def \horizendproof{\kern 5mm
\hbox{\vbox{\hrule\hbox to 5pt{\vrule height 4,2 pt\hfil\vrule}\hrule}
}}

\belowlistskipamount 12pt
\interitemskipamount 10pt

\titre Introduction

\definexref{quadratic functional}{\the\notitre}{titre}

In
signal processing, blind source separation consists in the
identification
of analogical, independent signals mixed by a black-box device.
In psychometrics, one has the notion of structural latent variable 
whose mixed effects are only measurable through series of tests; an example 
are the Big Five 
identified from factorial analysis 
by researchers in the domain of personality evaluation (Roch, 1995). 
Other application fields such as digital imaging, 
bio medicine, finance and econometrics also use  models aiming to 
recover hidden independent factors from observation. 
Independent component analysis (ICA) is one such tool; it can be seen as an extension
of principal component analysis, in that it goes beyond a simple 
linear decorrelation only satisfactory for a normal distribution; 
or as a complement, since its application is precisely pointless under the assumption of normality.

Papers on ICA are found in 
the  fields of signal processing,  neural networks, 
statistics and information theory. Comon (1994)
defined the concept of ICA as maximizing the degree of statistical independence among 
 outputs using contrast functions approximated by the Edgeworth expansion of 
 the Kullback-Leibler divergence.

The model is usually stated as follows: 
let  $X$ be a random variable on $\R^d$, $d \geq2$;   
find pairs $(A,S)$, such that $X=AS$, 
where $A$ is a square invertible matrix and $S$ 
a latent random variable whose components are mutually independent.
This is usually done by minimizing some contrast function that cancels
out if, and only if,
the components of $WX$ are independent, where $W$ is a candidate for the 
inversion of $A$.

Matrix $A$ is identifiable up to a scaling matrix and a  permutation matric   
if and only if $S$ has at most one Gaussian component (Comon, 1994).


Maximum-likelihood methods and contrast functions based on
mutual information or other divergence measures between densities are
commonly employed.
Bell and Snejowski (1990s) published an approach based on the Infomax principle. 
Cardoso (1999) used higher-order cumulant tensors, which led to the Jade algorithm,
Miller and Fisher III (2003) proposed a contrast based on a spacing estimates of entropy. 
Bach and Jordan (2002) proposed a contrast function based on canonical correlations in
a reproducing kernel Hilbert space.
 Similarly, Gretton et al (2003)
 proposed kernel covariance and kernel mutual information
contrast functions. 
Tsybakov and Samarov (2002)  proposed a method of direct estimation of $A$,  based on nonparametric estimates of  
matrix functionals using the gradient of $f_A$.

Let $f$ be the density of the latent variable $S$ relative to Lebesgue measure, assuming it exists. The observed variable 
$X=AS$ has the density $f_A$,
given by 
$$
\eqalign{
f_A(x)&=|\det A^{-1}| f(A^{-1}x)\cr
&= |\det B| f^1(b_1x)\ldots f^d(b_dx),\cr
}$$
where $b_\ell$ is the $\ell$th row  of the matrix $B=A^{-1}$; 
this resulting from a change of variable if the latent density
$f$ is equal to the product of its marginals $f^1\ldots f^d$. 
In this regard, 
 latent variable $S=(S^1\!,\ldots,S^d)$ having independent components means independence of 
the random variables $S^\ell\circ \pi^\ell$ defined
on some product probability space $\Omega=\prod \Omega^\ell$, with $\pi^\ell$ the canonical projections.
So $S$ can be defined as the compound of the
unrelated $S^1,\!\ldots,S^d$ sources.  

In the ICA model expressed this way, both $f$ and $A$ are unknown, and the data consists in a random sample of $f_A$. The semi-parametric case corresponds to $f$ left unspecified, except for general regularity assumptions.

\bigskip

In this paper, we consider the exact contrast provided by the   factorization 
measure $\int |f_A-f_A^\star|^2\,$,
with $f_A^\star$ the product of the marginals of $f_A$. 
Let's mention that the idea of comparing in the $L_2$ norm a joint density
with the product of its marginals, can be traced
back to Rosenblatt (1975).

\subsection Estimation of a quadratic functional

The problem of estimating nonlinear functionals of a density has been widely studied. 
In estimating $\int f^2$ under H\"older smoothness conditions, Bickel and Ritov (1988) have shown
 that  parametric rate  
is achievable for a regularity $s\geq 1/4$, whereas when $s\leq 1/4$,
minimax rates of convergence under mean squared error are of the order of $n^{-8s/1+4s}$. 
This result has been extended to general functionals of a density 
$\int \phi(f)$ by Birg\'e and Massart (1995). Laurent (1996) has built
efficient estimates for $s>1/4$.

Let $P_j$ be the projection operator on a multiresolution analysis (MRA) at
level $j$,
with  scaling function $\varphi$, and let $\alpha_{jk}=\int f\varphi_{jk}$ be the coordinate $k$ of $f$.

 In the wavelet setting, given a sample $\tilde X=\{X_1,\ldots, X_n\}$ of a density
$f$ defined on $\R$, independent and identically distributed, the U-statistic 
$$\hat B^2_j(\tilde X)={2\over n(n-1)}\sum_{i_1<i_2} \sum_{k\in\smallZ} \varphi_{jk}(X_{i_1})\varphi_{jk}(X_{i_2})$$
with mean $\int (P_jf)^2$ is the usual optimal estimator of the quantity
$\int f^2$ ; see Kerkyacharian and Picard (1996), and 
 Tribouley (2000) for the white noise model with adaptive rules. 

In what follows, this result is implicitly extended to $d$ dimensions using a tensorial wavelet basis $\Phi_{jk}$,
with $\Phi_{jk}(x)=\varphi_{jk^1}(x^1)\ldots\varphi_{jk^d}(x^d)$,
$k\in \Z^d, x\in \R^d$;
that is to say with $\tilde X$ an independent, identically distributed  sample of a density $f$ on $\R^d$, the
U-statistic
$\hat B^2_j(\tilde X)={2\over n(n-1)}\sum_{i_1<i_2} \sum_{k\in\smallZ^d} \Phi_{jk}(X_{i_1})\Phi_{jk}(X_{i_2})$
with mean $\int (P_jf)^2=\sum_{k\in \smallZ^d} \alpha_{jk}^2$ is also optimal in estimating 
the quantity $\int_{\smallR^d} f^2$.

In the case of a compactly supported density $f$, $\hat B^2_j$  is computable with a Daubechies wavelet D2N
and dyadic approximation of $X$, but  the computational cost is basically  in $O(n^2(2N-1)^d)$,  which is generally too  high in practice.

On the other hand, 
the plug-in,   biased, estimator 
$\hat H^2_j(f)=\sum_k \left[{1\over n}\sum \Phi_{jk}(X_i)\right]^2=\sum_k\hat \alpha_{jk}^2$ enjoys both 
ease of computation and ease of transitions between resolutions through discrete wavelet transform (DWT),
since it 
 builds upon a preliminary
estimation of all individual wavelet coordinates of $f$
on the projection space at level $j$, that is to say 
a full density estimation. 
In this setting it is just as easy to compute 
$\sum_k |\hat\alpha_{jk}|^p$ for any $p\geq 1$ or even $\sup |\hat\alpha_{jk}|$,
with a fixed computational cost in $O(n(2N-1)^d)$ 
plus sum total,  or
 seek out the max, of a $2^{jd}$ array.

Both estimators  $\hat H^2_j$ and $\hat B^2_j$ build
on the same  kernel 
$h_j(x,y)=\sum_{k\in\smallZ^d} \Phi_{jk}(x)\Phi_{jk}(y)$ since they  are written
$$
\hat H^2_j(\tilde X)= (n^2)^{-1} \sum_{i \in \Omega_n^2} h_j(X_{i^1},X_{i^2}) \quad \hbox{ and }\quad
\hat B^2_j(\tilde X)= (A_n^2)^{-1} \sum_{i \in I_n^2} h_j(X_{i^1},X_{i^2}),$$
where, here and in the sequel, 
$\Omega_n^m=\Onm$, $I_n^m=\InminOnm$ and
$A_n^p=n!/(n-p)!$.

The plug-in estimator $\hat H^2_j$ is then  identified as the Von Mises statistic associated to $\hat B^2_j$. In estimating $\sum_k \alpha_{jk}^2$,  the mean squared error of unbiased $\hat B^2_j$ is merely its variance, while the  mean squared error of $\hat H^2_j$ adds
a squared component $E (\hat H^2_j - \hat B^2_j)^2$
because of the inequality 
$(\hat H^2_j - \sum_k \alpha_{jk}^2)^2 \leq  2(\hat H^2_j - \hat B^2_j)^2 + 2(\hat B^2_j - \sum_k \alpha_{jk}^2)^2 $.

From general results, 
a U-statistic with finite second raw moment
has a variance in $Cn^{-1}$ and under similar conditions, the difference
$E | U-V|^r$ between the U-statistic and its associated Von Mises statistic is of the order of $n^{-r}$ (See for instance Serfling, 1980).

In the wavelet case, the dependence of the statistics on the resolution $j$ calls for special treatment in computing these two quantities. This special computation, taking $j$ and other properties of wavelets into account, constitutes 
the main topic of the paper. In particular whether
 $2^{jd}$ is lower than $n$ or not is a critical  threshold for resolution parameter $j$.
 Moreover, on the set $\{\,j\colon2^{jd} > n^2\}$, the statistic $\hat B^2_j$,  and therefore also $\hat H^2_j$, have a mean squared error not converging to zero.

If $\hat B^2_j$ and $\hat H^2_j$ share some features
in estimating $\sum_k \alpha_{jk}^2=\int (P_jf)^2$, they 
differ in an essential way: the kernel $h_j$ is averaged in one case over $\Omega_n^2$, the set of 
unconstrained indexes, and in the other case over $I_n^2$ the set of
distinct indexes.
As a consequence, it is shown in the sequel that $\hat H^2_j$  has mean squared error of the order of $2^{jd}n^{-1}$, which makes it inoperable
 as soon as  $2^{jd}\geq n$, while $\hat B^2_j$ 
 has mean squared error of the order of $2^{jd}n^{-2}$, which is then
 parametric on the set $\{\,j\colon 2^{jd}<n\}$.
 In a general way, this same parallel $\Omega_n^m$ versus $I_n^m$ is 
underpinning most of the proofs presented throughout the paper.

\subsection  Wavelet ICA 

Let $f$ be the latent density in the semi-parametric model introduced above. Let $f_A$ be the mixed density
and let $f_A^\star$ be the product of the marginals of $f_A$.

Assume, as regularity condition, that $f$ belongs to a Besov class $B_{s2\infty}$. It has been checked in previous work (Barbedor, 2005)
that $f_A$ and $f_A^\star$, hence $f_A-f_A^\star$ belong to the same Besov space than $f$.

As usual, the very definition of Besov spaces (here $B_{s2\infty}$) and an orthogonality property of the projection spaces $V_j$ and $W_j$ 
entails the relation 
$$0 \leq \int (f_A-f_A^\star)^2 - \int \left[P_j(f_A-f_A^\star)\right]^2 \leq C 2^{-2js}.$$
In this relation, the quantity $\int [P_j(f_A-f_A^\star)]^2$ is recognized as the wavelet ICA contrast $C^2_j(f_A-f_A^\star)$, introduced in a preliminary paper (Barbedor, 2005).

The wavelet ICA contrast is then a factorization measure with 
 bias,
in the sense that a zero contrast implies independence of 
the projected densities, and that  independence in projection
transfers to original densities up to some bias $2^{-2js}$.

Assume for a moment that the difference $f_A-f_A^\star$ is a density
and that we dispose of an independent, identically distributed sample $\tilde S$ of this difference.
Computing the estimators $\hat B^2_j(\tilde S)$ or $\hat H^2_j(\tilde S)$  
provides an estimation of $\int (f_A-f_A^\star)^2$, the exact ICA factorization measure. In this case,
the $j^\ast$ realizing the best compromise between  the mean squared
error in $C^2_j$ estimation and the bias of the ICA wavelet contrast $2^{-2js}$,
 is exactly the same as the one to minimize the overall risk in estimating 
the quadratic functional $\int (f_A-f_A^\star)^2$.
It is found by balancing bias and variance, a standard procedure in nonparametric estimation.
 From what was said 
above $\hat B^2_j(\tilde S)$ would be an optimal estimator of the exact factorization measure $\int (f_A-f_A^\star)^2$.

\bigskip The previous assumption being heuristic only, and since, in  ICA,  the data at hand is a random sample of $f_A$ and not $f_A-f_A^\star$, we are lead to consider estimators different from $\hat B^2_j$ and $\hat H^2_j$, but still alike in some way.

Indeed, let $\delta_{jk}=\int (f_A-f_A^\star)\Phi_{jk}$ be the coordinate of the difference function $f_A-f_A^\star$. In the ICA context, $\delta_{jk}$ is estimable only through the difference
$(\alpha_{jk} - \prodalpha)$ where $\alpha_{jk}=\int f_A\Phi_{jk}$
is the coordinate of $f_A$ and $\alpha_{jk^\ell}=\int\fix(f_A,\ell)\varphi_{jk^\ell}$
refers to the coordinate of  marginal number $\ell$ of $f_A$, written $\fix(f_A,\ell)$. 

To estimate $\sum_k \delta_{jk}^2$, estimators of the type $\hat B^2_j$
and  $\hat H^2_j$ are not alone enough. Instead
we use the already introduced wavelet contrast estimator (plug-in),
$\hat C^2_j(\tilde X)=\sum_{k} (\hat\alpha_{jk^1,\ldots k^d}
-\hat\alpha_{jk^1}\ldots \hat\alpha_{jk^d})^2$, and
the corresponding U-statistic estimator of order $2d+2$,
$$\eqalign{
\hat D^2_j(\tilde X)={1\over A_n^{2d+2}}
\sum_{i\in I_n^{2d+2}}\sum_{k\in\smallZ^d}  
\big[\Phi_{jk}(X_{i^1})&-\varphi_{jk^1}(X_{i^2}^1)\ldots
\varphi_{jk^d}(X_{i^{d+1}}^d)\big]\cr
&\big[\Phi_{jk}(X_{i^{d+2}})-
\varphi_{jk^1}(X_{i^{d+3}}^1)\ldots
\varphi_{jk^d}(X_{i^{2d+2}}^d)\big]
}$$ with as above $I_n^m=\Inm$ and $X^\ell$ referring to the dimension $\ell$ of $X\in \R^d$.

As it turns out, the U-statistic estimator $\hat D^2_j$ computed
on the full sample $\tilde X$ is slightly suboptimal, compared to the rate of a $\hat B^2_j$ in estimating  a bare quadratic functional.

As an alternative to $\hat D^2_j(\tilde X)$, we are then led to consider
 various U-statistic and plug-in estimators 
based on splits of the full sample, which seems the only way to find back the well-known optimal convergence rate of the estimation of quadratic functional, for reasons that will be explained in the course of the proofs. 

These additional estimators and conditions of use, together with the 
full sample estimators $\hat C^2_j$ and $\hat D^2_j$ are presented in section 3.

Section 2 of the paper recalls some essential definitions  for the convenience of the reader not familiar with wavelets and Besov spaces, and may be skipped.

Section 4 is all devoted to the computation of a risk bound
for the different estimators  presented in section 3.

We refer the reader to a preliminary paper on ICA by wavelets (Barbedor, 2005) which contains numerical simulations, details
on the implementation of the wavelet contrast estimator and other practical considerations not repeated here. Note that this paper gives an improved convergence rate in $C2^{jd}n^{-1}$ for the wavelet contrast estimator $\hat C^2_j$, already introduced in the preliminary paper.

\section Notations

We set here general notations and recall some definitions for  the convenience of ICA specialists. The reader already familiar with wavelets and Besov spaces can skip this part.

\unorderedlist

\li {\it Wavelets} 

Let $\varphi$ be some function of $L_2(\R)$ such that
the family of translates $\{\varphi(.-k),\, k\in \Z\}$ is an orthonormal 
system; let $V_j\subset L_2(\R)$ be the subspace spanned by 
 $\{\varphi_{jk}= 2^{j/2}\varphi(2^j . -k), k \in \Z\}$.

By definition, the sequence of spaces $(V_j), j\in \Z$,
is called a multiresolution analysis (MRA) of $L_2(\R)$ if
$V_j \subset V_{j+1}$ and $\bigcup_{j\geq 0} V_j$ is dense in $L_2(\R)$;
$\varphi$ is called the father wavelet or scaling function.

Let $(V_j)_{j\in \smallZ}$ be a multiresolution analysis
of $L_2(\R)$, with $V_j$ spanned by 
 $\{\varphi_{jk}= 2^{j/2}\varphi(2^j . -k), k \in \Z\}$. 
Define $W_j$  as the complement of  $V_j$ in $V_{j+1}$, and let the families
$\{\psi_{jk}, k\in \Z\}$ be a basis for $W_j$,
with $\psi_{jk}(x)= 2^{j/2}\psi(2^j x -k)$. 
Let  $\alpha_{jk}(f)=<f,\varphi_{jk}>$ and $\beta_{jk}(f)=<f,\psi_{jk}>$.
 
A function $f\in L_2(\R)$ admits a wavelet expansion on $(V_j)_{j\in \smallZ}$ if
 the series
$$\sum\limits_{k} \alpha_{j_0k}(f) \varphi_{jk}
+ \sum\limits_{j=j_0}^\infty
\sum\limits_{k} \beta_{jk}(f) \psi_{jk}
$$ is convergent to $f$ in $L_2(\R)$;
$\psi$ is called a mother wavelet.
\bigskip

A MRA in dimension one also induces an associated 
MRA in dimension $d$, using the tensorial product procedure below.

Define $V^d_j$ as the tensorial product
of $d$ copies of $V_j$.
The increasing 
sequence $(V^d_j)_{j\in \smallZ}$
defines a multiresolution analysis of $L_2(\R^d)$
(Meyer, 1997):

-- for $(i^1\ldots,i^d) \in \{0,1\}^d$
and $(i^1\ldots,i^d) \neq (0\ldots,0)$,  
define 
$$\Psi(x)_{i^1\ldots,i^d}=\prod_{\ell=1}^d \psi^{(i^\ell)}(x^\ell),
\eqdef{psix}$$
 with $\psi^{(0)}=\varphi$, $\psi^{(1)}=\psi$,
so that $\psi$ appears at least once in the product $\Psi(x)$
(we now on omit $i^1\ldots,i^d$ in the notation for $\Psi$, and
in \eqref{wavexpansion}, although it
is present each time);

-- for $(i^1\ldots,i^d) = (0\ldots,0)$,  define 
  $\Phi(x)=\prod_{\ell=1}^d \varphi(x^\ell)$;

-- for $j \in \Z$, $k\in \Z^d$, $x \in \R^d$, let  
$\Psi_{jk}(x)=2^{jd \over 2} \Psi(2^jx-k)$ and 
$\Phi_{jk}(x)=2^{jd \over 2} \Phi(2^jx-k)$;

-- define $W^d_j$ as the orthogonal complement of  
$V^d_j$ in $V^d_{j+1}$; 
it is an orthogonal sum of $2^d-1$ spaces
having the form $U_{1j}\ldots\otimes U_{dj}$, where $U$
is a placeholder for $V$ or $W$; $V$ or $W$ 
are thus placed using up all permutations, but with 
  $W$ represented at least once, so that 
a fraction of the overall innovation brought by
 the finer resolution $j+1$ is always present in the tensorial product.

A function $f$ admits a wavelet expansion on the basis $(\Phi, \Psi)$ 
if the series
$$\sum\limits_{k\in \smallZ^d} \alpha_{j_0k}(f) \Phi_{j_0k}
+ \sum\limits_{j=j_0}^\infty \sum\limits_{k\in \smallZ^d} \beta_{jk}(f) 
\Psi_{jk}
\eqdef{wavexpansion}$$ is convergent to $f$ in $L_2(\R^d)$.

\bigskip

In connection with function approximation, wavelets can be viewed 
as falling in the category of orthogonal
series methods, or also in the category of kernel methods.

The approximation at level $j$ of a function $f$ that admits a 
multiresolution
expansion is the orthogonal projection  $P_jf$ of $f$ onto
$V_j\subset L_2(\R^d)$  defined by
$$(P_jf)(x)=\sum_{k \in \smallZ^d}\alpha_{jk}\Phi_{jk}(x),
$$
where  
 $\alpha_{jk}=\mulalpha=\int f(x)\Phi_{jk}(x)\,dx$.

With a concentration condition verified for compactly supported wavelets, 
the projection operator can also be  written
$$
(P_jf)(x)=\int_{\smallR^d} K_j(x,y)f(y)d(y),
$$
with $K_j(x,y)=2^{jd}\sum_{k\in \smallZ^d}\Phi_{jk}(x)
\overline{\Phi_{jk}(y)}$. $K_j$ is an orthogonal projection kernel
 with window
$2^{-jd}$ (which is not translation invariant).

\li {\it Besov spaces} 

Besov spaces  admit a characterization in terms of wavelet coefficients,
which makes them intrinsically connected
to the analysis of curves via wavelet techniques.

 $f\in L_p(\R^d)$ belongs to the (inhomogeneous) Besov space $B_{spq}(\R^d)$ if
$$J_{spq}(f)=\| \alpha_{0.}\|_{\ell_p} + 
\bigg[ \sum_{j\geq 0}\left[  
2^{js}2^{dj({1\over 2}-{1\over p}})\|\beta_{j.} \|_{\ell_p}
\right]^q\bigg]^{1\over q} < \infty,$$
with  $s>0$, $1\leq p \leq\infty$, $1\leq q \leq\infty$, 
and $\varphi, \psi$ $\in C^r, r>s$ (Meyer, 1997).

Let $P_j$ be the projection operator on $V_j$ and let $D_j$ be the projection operator on $W_j$. 
$J_{spq}$ is equivalent to 
$$J'_{spq}(f)=\| P_jf\|_p + 
\bigg[ \sum_{j\geq 0}\left[  
2^{js}\|D_j f \|_p
\right]^q\bigg]^{1\over q}$$

\endunorderedlist
 
A more complete presentation of wavelets linked with Sobolev and Besov approximation theorems 
and statistical applications can be found in 
the book from H\"ardle et al. (1998). 
General references about Besov spaces
are  Peetre (1975), 
Bergh \& L\"ofstr\"om (1976), Triebel (1992), DeVore \& Lorentz (1993).

\section Estimating the factorization measure $\int (f_A-f_A^\star)^2$ 
 
We first recall the definition of the wavelet contrast already introduced in Barbedor(2005). 

Let $f$ and $g$ be two functions on $\R^d$ and let 
$\Phi$ be the
scaling function of a multiresolution analysis
of $L_2(\R^d)$ for which projections of $f$ and $g$ exist.

Define the approximate loss function 
$$C^2_j(f-g)=\sum_{k \in\smallZ^d}
\left(\int (f-g)\Phi_{jk}\right)^2=\|P_j(f-g)\|_2^2.$$
It is clear that $f=g$ implies  $C^2_j=0$ and that 
$C^2_j=0$ implies $P_jf=P_jg$ almost surely.

Let $f$ be a density function on $\R^d$; denote 
by $\fix(f,\ell)$ the marginal 
 distribution in dimension $\ell$
$$ x^\ell \mapsto \int_{\smallR^{d-1}}
f(\echd(x))\,dx^1\ldots dx^{\ell-1}dx^{\ell+1}\ldots dx^d$$
and denote by $f^\star$ the product of marginals
$\fix(f,1)\ldots\fix(f,d)$.
The functions $f$, $f^\star$ and the $\fix(f,\ell)$
 admit a wavelet expansion 
on a  compactly supported basis $(\varphi,\psi)$.
Consider the projections up to order $j$,  
that is to say the projections of $f$, $f^\star$ and $\fix(f,\ell)$ 
on  $V^{d}_j$ and $V_j$, namely
$$
P_jf^\star=\sum_{k \in \smallZ^d}\alpha_{jk}(f^\star)
\Phi_{jk},
 \quad  
P_jf=\sum_{k \in \smallZ^d}\alpha_{jk}(f)\Phi_{jk}
 \quad \hbox{ and } 
P_j^\ell\fix(f,\ell)=\sum_{k \in \smallZ}
\alpha_{jk}(\fix(f,\ell))\varphi_{jk},
$$
with
$\alpha_{jk}(\fix(f,\ell))
=\int \fix(f,\ell)\varphi_{jk}$
and $\alpha_{jk}(f)=\int f\Phi_{jk}$. 
At least for compactly supported densities and compactly supported wavelets, it is clear that  
$P_jf^\star=P_j^1\fix(f,1)
\ldots P_j^d\fix(f,d)$.

 


\proposition ICA wavelet contrast\par
\definexref{wavelet contrast}{\the\notitre.\the\noproposition}{proposition}
{
\it
\noindent
Let $f$ be a compactly supported density function on $\R^d$ and let
 $\varphi$ be the  
scaling function of a compactly supported wavelet.

Define the wavelet ICA contrast as $C^2_j(f-f^\star)$.
Then,  
$$\eqalign{
 f \hbox{ factorizes} &\kern 7pt\Longrightarrow \kern 5pt C^2_j(f-f^\star)=0\cr
  C^2_j(f-f^\star)=0 &\kern 7pt\Longrightarrow \kern 5pt P_jf = P_j\fix(f,1)\ldots P_j\fix(f,d)\quad a.s.
}$$

}

{\bf Proof}
  $f = f^1\ldots f^d \Longrightarrow \fix(f,\ell)=f^\ell$, $\ell=1,\ldots d$. 
\endproof


\subsection Wavelet contrast and quadratic functional

Let $f=f_I$ be a density defined on $\R^d$ whose components are independent, that is to say
$f$ is equal to the product of its marginals.
Let $f_A$ be the mixed density 
given by 
$
f_A(x)=|\det A^{-1}| f(A^{-1}x)
$, with $A$ a $d\times d$ invertible matrix. Let $f_A^\star$ be the product of the marginals of $f_A$. Note that when $A=I$, $f_A^\star=f_I^\star=f_I=f$.

By definition of a Besov space $B_{spq}(\R^d)$ 
with a $r$-regular wavelet $\varphi$, $r > s$,
$$ f\in B_{spq}(\R^d) \Longleftrightarrow
\| f-   P_j f\|_p = 2^{-js}\, \epsilon_j, \quad \{\epsilon_j\} \in \ell_q(\N^d).\eqdef{meyer}
$$

So, from the decomposition
$$\eqalign{
\| f_A  - f_A^\star\|_2^2 &= \int P_j(f_A  - f_A^\star)^2 + \int 
\big[ f_A- f_A^\star - P_j(f_A- f_A^\star)\big]^2,
\cr
&= C^2_j(f_A  - f_A^\star) + \int 
\big[ f_A- f_A^\star - P_j(f_A- f_A^\star)\big]^2,
\cr
}$$ resulting from the  orthogonality of $V_j$ and $W_j$,
and assuming that $f_A$ and $f_A^\star$ belong to $B_{s2\infty}(\R^d)$,
$$0 \leq \| f_A  - f_A^\star\|_2^2 - C^2_j(f_A-f_A^\star) \leq C 2^{-2js},
\eqdef{ketoile}$$

which gives an illustration of the shrinking (with $j$) 
distance between 
the wavelet contrast and the always bigger squared $L_2$ norm of 
$f_A-f_A^\star$ representing the exact factorization measure.  
A side effect of \eqref{ketoile} is
that $C^2_j(f_A-f_A^\star)=0$ is implied by $A=I$.



\subsection  Estimators under consideration

Let $S$ be the latent random variable with density $f$.\br
Define the  experiment
${\cal E}^n=\bigl({\cal X}^{\otimes n},\, {\cal A}^{\otimes n},\,
(X_1, \ldots, X_n),\,
\Pnfa,\, f_A \in B_{spq}\bigr)$, where $X_1, \ldots, X_n$ is an iid sample of $X=AS$, 
and $\Pnfa=\Pfa\ldots \otimes \Pfa$
is the joint distribution of ($\echn(X)$).

Define the coordinates estimators
$$\hat\alpha_{jk}=\hat \mulalpha = 
{1 \over n} \sum\limits_{i=1}^n \varphi_{jk^1}(X_i^1)\ldots\varphi_{jk^d}(X_i^d)
\quad\hbox{ and }\quad\hat\alpha_{jk^\ell}= 
{1 \over n} \sum\limits_{i=1}^n \varphi_{jk^\ell}(X_i^\ell)
\eqdef{defcoordinates}$$
where $X^\ell$ is coordinate $\ell$ of $X\in \R^d$. Define also the shortcut
$\hat\lambda_{jk}=\hat\alpha_{jk^1}\ldots\hat\alpha_{jk^d}$.

Define the full sample plug-in estimator 
$$\hat C^2_j=\hat C^2_j(X_1,\ldots,X_n)=\sum_{(k^1,\ldots,k^d)\in \smallZ^d}(\hat\alpha_{j(k^1,\ldots,k^d)}-\hatprodalpha)
^2
=\sum_{k\in \smallZ^d}(\hat\alpha_{jk}-\hat\lambda_{jk})^2\eqdef{defC2j}$$ and the full sample U-statistic estimator
$$\eqalign{
\hat D^2_j=\hat D^2_j(X_1,\ldots,X_n)={1\over A_n^{2d+2}}
\sum_{i\in I_n^{2d+2}}\sum_{k\in\smallZ^d}  
\big[\Phi_{jk}(X_{i^1})&-\varphi_{jk^1}(X_{i^2}^1)\ldots
\varphi_{jk^d}(X_{i^{d+1}}^d)\big]\cr
&\big[\Phi_{jk}(X_{i^{d+2}})-
\varphi_{jk^1}(X_{i^{d+3}}^1)\ldots
\varphi_{jk^d}(X_{i^{2d+2}}^d)\big]
}\eqdef{defD2j}$$
where $I_n^m$ is the set of indices $\Inm$ and $A_n^m={n!\over(n-m)!}=|I_n^m|$.
%

Define also the U-statistic estimators 
$$\eqalign{\hat B^2_j(\{X_1,\ldots,X_n\})&=\sum_{k} {1\over A_n^2}\sum_{i\in I_n^2}\Phi_{jk}(X_{i^1})\Phi_{jk}(X_{i^2}) \cr
\hat B^2_j(\{X_1^\ell,\ldots,X_n^\ell\})&=\sum_{k^\ell} {1\over A_n^2}\sum_{i\in I_n^2}\varphi_{jk^\ell}(X^\ell_{i^1})\varphi_{jk^\ell}(X^\ell_{i^2}).}
\eqdef{B2Ustat}$$

\subsection  Notational remark

 Unless
otherwise stated, superscripts
designate coordinates of multi-dimensional entities while
subscripts designate unrelated entities of the same set without reference to multi-dimensional unpacking.
For instance, an index $k$ belonging to $\Z^d$ is also written  $k=(k^1,\ldots,k^d)$, with $k^\ell\in \Z$. Likewise
 a multi-index $i$ is written $i=(i^1,\ldots,i^m)$ when belonging to some $\Omega_n^m=\{i=(i^1,\ldots,i^m)\colon i^\ell\in \N,\, 1\leq i^\ell\leq n\}$
 or $I_n^m=\{i\in \Omega_n^m\colon \ell_1\neq \ell_2\Rightarrow i^{\ell_1}\neq i^{\ell_2} \}$, for some $m\geq 1$; but $i_1$, $i_2$ would designate two different elements of $I_n^m$,
so  for instance 
$\left[\sum_{i=1}^n\sum_{k\in \smallZ^d} \Phi_{jk}(X_i)\right]^2$ is written  $\sum_{ i_1,i_2}\sum_{k_1,k_2}\Phi_{jk_1}(X_{i_1})\Phi_{jk_2}(X_{i_2})$. Finally  $X^\ell$ is coordinate $\ell$ of observation $X\in \R^d$ and $\tilde X$ refers to a sample $\{X_1,\ldots X_n\}$.\puncrem

As was said in the introduction and as is shown in proposition
\ref{onesampleD2j}, the estimator $\hat D^2_j$ computed on the full sample is slightly suboptimal.
We now  review some  possibilities to split the sample so that various
alternatives to $\hat D^2_j$ on the full sample could be computed
in an attempt to regain optimality through block independence. 

We need not consider $\hat C^2_j$ on independent sub samples
because, as will be seen, the order of its risk upper bound is given by the
order of the component $\sum_k \hat\alpha_{jk}^2-\alpha_{jk}^2$ which is not improved by splitting the sample (contrary to 
$\sum_k \hat\lambda_{jk}^2-\lambda_{jk}^2$ and
$\sum_k \hat\alpha_{jk}\hat\lambda_{jk}-\alpha_{jk}\lambda_{jk}$).
The rate of $\hat C_j^2$ is unchanged compared to
what appeared in Barbedor (2005).

\subsection Sample split

\medskip
\unorderedlist
\li
Split the full sample $\{X_1,\ldots, X_n\}$ in $d+1$ disjoint sub samples $\tilde R^0, \tilde R^1,\ldots \tilde R^d$
where the sample $\tilde R^0$ refers to a plain section of the full sample,
$\{X_1,\ldots, X_{[n/d+1]}\}$ say,
and the samples $\tilde R^1,\ldots,\tilde R^d$ refer to dimension $\ell$ of their section of the full sample, for instance $\{X^\ell_{[n/d+1]\ell+1},\ldots, X^\ell_{[n/d+1](\ell+1)}\}$.

Estimate each plug-in $\hat \alpha_{jk}(\tilde R^0)$ and $\hat \alpha_{jk^\ell}(\tilde R^\ell)$, and the U-statistics $\hat B^2_j(\tilde R^0)$, $\hat B^2_j(\tilde R^\ell)$, $\ell=1,\ldots,d$ on each independent sub-sample. This leads to the definition of the $d+1$ samples 
mixed plug-in estimator
$$
\hat F^2_j(\tilde R^0,\tilde R^1,\ldots,\tilde R^d)=\hat B^2_j(\tilde R^0) 
+ \prod_{\ell=1}^d \hat B^2_j(\tilde R^\ell)-2\sum_{k\in\smallZ^d}\hat\alpha_{jk}(\tilde R^0)
\hat\alpha_{jk^1}(\tilde R^1)\ldots\hat\alpha_{jk^d}(\tilde R^d).
\eqdef{mixedplugind+1sample}$$
to estimate the quantity
$\sum_k\alpha_{jk}^2 + \prod_{\ell=1}^d\Bigl(\sum_{k^\ell\in\smallZ}\alpha_{jk^\ell}^2\Bigr)
-2 \sum_k\alpha_{jk}\alpha_{jk^1}\ldots\alpha_{jk^d}=C^2_j$.

Using estimators $\hat B^2_j$  places us 
in the exact replication of the case $\hat B^2_j$ 
found in Kerkyacharian and Picard (1996), except for an estimation taking place in dimension $d$ in the case of $\hat B^2_j(R^0)$.
The risk of this procedure is given by proposition 
\ref{d+1subsample}.

\li
Using the full sample $\{X_1,\ldots,X_n\}$ we can generate 
an identically distributed sample of $f_A^\star$,
namely $\tilde{DS}=\cup_{i\in \Omega_n^d} \{X_{i^1}^1\ldots X_{i^d}^d\}$,
but is not constituted of independent observations when $A\neq I$. 

But then using a Hoeffding like decomposition, we can pick from $\tilde{DS}$, 
a sample of independent observations,  
$\tilde{IS}=\cup_{k=1\ldots [n/d]} \{X_{(k-1)d+1}^1\ldots X_{kd}^d\}$,  although it leads to a somewhat arbitrary omission
of a large part of the information available.
Nevertheless  we can  assume that we dispose of two independent, identically distributed samples, one for $f_A$ labelled $\tilde R$ and one for $f_A^\star$ labelled $\tilde S$, with $\tilde R$ independent of $\tilde S$. In this setting
we define the mixed plug-in estimator
$$\hat G^2_j(\tilde R,\tilde S)=\hat B^2_j(\tilde R)+\hat B^2_j(\tilde S)-2\sum_{k\in\smallZ^d}\hat\alpha_{jk}(\tilde R)\hat\alpha_{jk}(\tilde S)
\eqdef{mixedplugin2sample}$$
 and the two samples U-statistic estimator
$$ 
\hat \Delta^2_j(\tilde R,\tilde S)={1\over A_n^{2}}
\sum_{i\in I_n^{2}}\sum_{k\in \smallZ^d}  \big[\Phi_{jk}(R_{i^1})-\Phi_{jk}(S_{i^1})\big]
 \big[\Phi_{jk}(R_{i^2})-\Phi_{jk}(S_{i^2})\big]
\eqdef{Ustat2sample}$$
assuming for simplification that both samples have same size $n$
(that would be different from the size of the original sample).
$\hat \Delta^2_j(R,S)$ is the exact replication (except for dimension $d$ instead of 1) of the optimal estimator of $\int (f-g)^2$ for unrelated $f$ and $g$ found in Butucea and Tribouley (2006). The risk of this optimal procedure is found in   proposition \ref{2subsample}.
\endunorderedlist
\punclist

\subsection Bias variance trade-off

Let an estimator $\hat T_j$
be used in estimating the quadratic functional 
$K_\star=\int (f_A-f_A^\star)^2$; 
using \eqref {ketoile}, an upper bound for the mean squared error 
of this procedure when $f_A \in B_{s2\infty}(\R^d)$ is given by
$$\eqalign{
\Enfa(\hat T_j - K_\star)^2 
&\leq 2 \Enfa (\hat T_j - C^2_j)^2 +
 C 2^{-4js},
}\eqdef{deuxjsbalance}$$
which shows that the key estimation is that of the wavelet 
contrast $C^2_j(f_A-f_A^\star)$
by the estimator $\hat T_j$. Once an upper bound of the risk of $\hat T_j$ in estimating $C^2_j$ is known, balancing the order of the bound with
the squared bias $2^{-4js}$ gives the optimal resolution $j$.
This is a standard procedure in nonparametric estimation.

Before diving into  the computation of risk bounds,  we give
a summary of the different convergence rates
in proposition \ref{linear-choice} below. The estimators based on splits
of the full sample are optimal. $\hat D^2_j$ is almost parametric
on $\{2^{jd} < n\}$ and is otherwise optimal.

\proposition Minimal risk resolution in the class $B_{s2\infty}$ and convergence rates\par
\definexref{linear-choice}{\the\notitre.\the\noproposition}{proposition}
{\it
Assume that $f$ belongs to $B_{s2\infty}(\R^d)$,
and  that projection is based on a $r$-regular wavelet $\varphi$, $r > s$.
Convergence rates for the estimators defined at the beginning of this section are the following:
$$\vbox{\smalltt
\offinterlineskip\hrule\halign {&\vrule#&\mystrut\quad\hfil#\quad\cr
depth 6pt&\omit&\omit&\omit&\omit&\omit&\cr
&\multispan5\hfil Convergence rates\hfil&\cr
depth 4pt&\omit&\omit&\omit&\omit&\omit&\cr
\noalign{\hrule}
depth 4pt&\omit&&\omit&&\omit&\cr
&statistic&&$2^{jd}< n$&&$2^{jd}\geq n$&\cr
height 2pt&\omit&&\omit&&\omit&\cr
\noalign{\hrule}
height 2pt depth 4pt&\omit&&\omit&&\omit&\cr
&$\hat \Delta_j^2(\tilde R,\tilde S),\; \hat G^2_j(\tilde R,\tilde S),\; \hat F^2_j(\tilde R^0,\tilde R^1,\ldots ,\tilde R^d)$&&parametric&&$n^{-8s\over 4s+d}$&\cr
&$\hat D_j^2(\tilde X)$&&$n^{-1+{1\over 1+4s}}$&&$n^{-8s\over 4s+d}$&\cr
&$\hat C_j^2(\tilde X)$&&$n^{{-4s\over 4s+d}}$&&inoperable&\cr
height 2pt&\omit&&\omit&&\omit&\cr
\noalign{\hrule}
}}  
$$
\centerline{\smallsf  Table 7. Convergence rates at optimal $j_\star$
}

The minimal risk resolution $j_\star$ satisfies,
 $2^{j_\star d}\approx(<)  n$ for parametric cases ; 
$2^{j_\star d} \approx n^{ 1 + {d - 4s \over d +4s }}$ for $\hat D^2_j$, $\hat \Delta^2_j$, $\hat G^2_j$ or $\hat F^2_j$ when $s\leq {d\over4}$ and
 $2^{j_\star d} \approx 
n^{ d \over d +4s }$ for $\hat C^2_j$.
}

\proof

Besov assumption about $f$ transfers to $f_A$ (see Barbedor, 2005).
Using 
$$\eqalign{
\Enfa(\hat H_j - K_\star)^2 
&\leq 2 \Enfa (\hat H_j - C^2_j)^2 +
 C 2^{-4js},
}$$
and balancing 
bias $2^{-4js}$ and variance of the estimator $\hat H_j$,
yields the optimal resolution $j$.

\unorderedlist

\li from proposition \ref{onesampleC2j}, for estimator $\hat C^2_j(\tilde X)$, the bound is inoperable
on $\{2^{jd}>n\}$. Otherwise equating $2^{jd}n^{-1}$ with
$2^{-4js}$ yields $2^j=n^{1\over d+4s}$ and a rate in 
$n^{{-4s\over d+4s}}$. 
\li from proposition \ref{2subsample} and \ref{d+1subsample}, for estimators $\hat F^2_j(R^0, R^1,\ldots,R^d)$, 
$\hat F^2_j(R,S)$ and $\hat D^2_j(R,S)$ , 
on $\{2^{jd} > n\}$ equating $2^{jd}n^{-2}$ with $2^{-4js}$ 
yields $2^j=n^{2\over d+4s}$ and a rate in 
$n^{-{8s\over d+4s}}$; on $\{2^{jd} < n\}$ the rate is  parametric.
Moreover $2^{jd}< n$ implies that $s\geq d/4$
and  $2^{jd}> n$ implies that $s\leq d/4$.
\li from proposition \ref{onesampleD2j}, for estimator 
$\hat D^2_j(\tilde X)$ 
on $\{2^{jd} > n\}$ equating $2^{jd}n^{-2}$ with $2^{-4js}$ 
yields $2^j=n^{2\over d+4s}$ and a rate in 
$n^{-{8s\over d+4s}}$; on $\{2^{jd} < n\}$ 
the rate is  found by equating 
$2^jn^{-1}$ with $2^{-4js}$. 

\endunorderedlist
\endproof
\section Risk upper bounds in estimating the wavelet contrast

In the forthcoming lines, we  make the assumption that both the density and the wavelet
are compactly supported so that all sums in $k$ are finite. For
simplicity we further suppose the density support to be the hypercube, so that
$\sum_{k\in\smallZ^d} 1 \approx 2^{jd}$.

\proposition  Risk upper bound, $d+1$ independent samples --- 
$f_A,\, \fix(f_A, 1),\ldots, \fix(f_A,d)$\par
\definexref{d+1subsample}{\the\notitre.\the\noproposition}{proposition}

{\it
Let $\{X_1,\ldots, X_n\}$ be an independent, identically distributed sample of $f_A$.
Let $\{R^\ell_1,\ldots, R^\ell_n\}$ be  an independent, identically distributed sample of $\fix(f_A,\ell)$, $\ell=1,\ldots,d$.  Assume that $f$ is compactly supported
 and that $\varphi$ is a Daubechies  $D2N$.
Assume that the $d+1$ samples are independent.
Let $\Enfa$ be  the expectation relative to the joint distribution
of the d+1 samples.
Then 
on $\{2^{jd}<n^2\}$,
 $$ \eqalign{
\Enfa\left(\hat F^2_j(\tilde X,\tilde R^1,\ldots,\tilde R^d) -C^2_j\right)^2 &\leq Cn^{-1}+C2^{jd} n^{-2}\indic{2^{jd}>n}. }$$ 
 
}
\proof

For the U-statistic $\hat F^2_j(\tilde X,\tilde R^1,\ldots,\tilde R^d)$, with $\hat\alpha_{jk}=\hat\alpha_{jk}(\tilde X)$, $\hat\alpha_{jk^\ell}=\hat\alpha_{jk^\ell}(\tilde R^\ell)$
and $\hat\lambda_{jk}=\hatprodalpha$,
$$\eqalign{
 (\hat F^2_j -C^2_j)^2&\leq 3\Bigl[\hat B^2_j(\tilde X)-\sum_k\alpha_{jk}^2\Bigr]^2
+ 3  \Bigl[\prod_\ell \hat B^2_j(\tilde R^\ell)
-\prod_\ell \sum_{k^\ell}\alpha_{jk^\ell}^2\Bigr]^2 + 
6\Bigl[\sum_k \hat\alpha_{jk}\hat\lambda_{jk}
-\sum_k \alpha_{jk}\lambda_{jk}\Bigr]^2.
}$$
On $\{2^{jd}<n^2\}$, by proposition \ref{variance hatb2j} for the term on the left, proposition  \ref{variance lambda deux} for the middle term, and proposition \ref{indepmultisample} for the term on the right,
the quantity is bounded by $Cn^{-1}+C2^{jd} n^{-2}$. 


\endproof

\proposition  Risk upper bound, 2 independent samples --- $f_A,\, f_A^\star$\par
\definexref{2subsample}{\the\notitre.\the\noproposition}{proposition}

{\it
Let $\tilde X=\{X_1,\ldots, X_n\}$ be an independent, identically distributed sample of $X$ with density $f_A$.
Let $\tilde R=\{R_1,\ldots, R_n\}$ be an independent, identically distributed sample of $R$ with density $f_A^\star$.  Assume that $f$ is compactly supported
 and that $\varphi$ is a Daubechies  $D2N$.
Assume that the two samples are independent.
Let $\Enfa$ be the the expectation relative to the joint distribution
of the two samples.

Then
 $$\eqalign{
   \Enfa\left(\hat G^2_j(\tilde X,\tilde R) -C^2_j\right)^2 &\leq Cn^{-1}+C2^{jd} n^{-2}\indic{2^{jd}>n}\cr
  \Enfa\left(\hat \Delta^2_j(\tilde X,\tilde R) -C^2_j\right)^2 &\leq C^\star n^{-1}+C2^{jd} n^{-2}. }$$with $C^\star=0$ at independence. 
 
}
\proof 

For the estimator
 $\hat G^2_j(\tilde X,\tilde R)$ the proof is identical to the proof of  proposition
 \ref{d+1subsample}, the only difference being that 
 $\hat\lambda_{jk}$ and $\lambda_{jk}$ no more designate a product of $d$ one dimensional coordinates but full fledged $d$ dimensional coordinate equivalent to $\hat\alpha_{jk}$ and $\alpha_{jk}$.
 
The only new quantity to compute is  then
$\Enfa\Bigl(\sum_k \hat\alpha_{jk}(\tilde X)\hat\lambda_{jk}(\tilde R)-\sum_k\alpha_{jk}\lambda_{jk}\Bigr)^2$, coming from the crossed term. 

Let $Q=\Enfa \bigl(\sum_k \hat\alpha_{jk}(\tilde X)\hat\lambda_{jk}(\tilde R)\bigr)^2$. Let $\theta=\sum_k \alpha_{jk}\lambda_{jk}$. 
Recall that $\Omega_n^m=\Onm$.

Let $\tilde \imath$ be the set of distinct coordinates of $i\in \Omega_n^4$.
So that, estimators being plug-in, with a sum on $\Omega_n^4$, with cardinality $n^4$, 
$$\eqalign{
Q
&=\Enfa{1\over n^4}\sum_{i\in\Omega_n^4}\sum_{k_1,k_2}\Phi_{jk_1}(X_{i^1})
\Phi_{jk_1}(R_{i^2})\Phi_{jk_2}(X_{i^3})\Phi_{jk_2}(R_{i^4})\cr
&\leq {1\over n^4}\bigg[\sum_{|\tilde\imath|=4} \theta^2+ 
\sum_{|\tilde\imath|=3} \Bigl[\theta^2 + 
(4N-3)^d\sum_k  \Enfa\Phi(X)^2\lambda_{jk}^2+
(4N-3)^d\sum_k  \Enfa\alpha_{jk}^2\Phi(R)^2\Bigr] + \cr
&\phantom{{1\over n^4}\bigg[}+\sum_{|\tilde\imath|\leq 2} 
(4N-3)^d\sum_k  \Enfa\Phi(X)^2\Phi(R)^2\bigg]\cr
}$$
with lines 2 and 3
expressing all possible matches between the coordinates of $i$,
and using lemma \ref{dbconcentration} to reduce double sums in $k_1,k_2$.

By independence of the samples, using lemma \ref{raw moments varphi}
and the fact that  $|\{i\in\Omega_n^4\colon |\,\tilde\imath\,|=c\}|=O(n^c)$ given by lemma \ref{generalized matching indices},
$$\eqalign{
Q
&\leq {A_n^4\over n^4}\theta^2+ 
Cn^{-1}  \Bigl(\theta^2 + C \sum_k \lambda_{jk}^2
+ C \sum_k \alpha_{jk}^2\Bigr)+ Cn^{-2} 2^{jd}.\cr
}$$
with $A_n^p=n!/(n-p)!$. So that, with $A_n^4n^{-4}=1 -{6\over n} + Cn^{-2}$,
$$\eqalign{
Q - \theta^2
&\leq 
Cn^{-2}  + Cn^{-1} + Cn^{-2} 2^{jd}.\cr
}$$

The rate is thus unchanged for $\hat F^2_j$ compared
to the $d+1$ sample case in previous proposition.

{\bf Case $\hat \Delta^2_j(\tilde X,\tilde R)$}

Recall that $I_n^m=\Inm$.

For $i\in I_n^2$, 
let
$h_{jk}(i)= \big[\Phi_{jk}(X_{i^1})-\Phi_{jk}(R_{i^1})\big]
\big[\Phi_{jk}(X_{i^2})-\Phi_{jk}(R_{i^2})\big]$ and let 
$\theta=C^2_j$;
so that 
$$\eqalign{
\Enfa \left(\hat \Delta^2_j(\tilde X, \tilde R) -\theta\right)^2
&= -\theta^2+ \Enfa{1\over (A_n^2)^2}\sum_{i_1,i_2} \sum_{k_1,k_2} h_{jk_1}(i_1)h_{jk_2}(i_2)
\cr
& =  \Bigl({\#\{i_1,i_2\colon|i_1\cap i_2|=0 \}
\over (A_n^2)^2}-1\Bigr) \theta^2 + 
{1\over (A_n^2)^2}\sum_{|i_1\cap i_2|\geq 1} \sum_{k_1,k_2}
\Enfa h_{jk_1}(i_1)h_{jk_2}(i_2)
,
\cr
 }$$
and by lemma \ref{matching indices} the quantity in parenthesis
on the left is of the order of $Cn^{-2}$.

Label $Q(i_1,i_2)$ the quantity $\Enfa  \sum_{k_1,k_2} h_{jk_1}(i_1)h_{jk_2}(i_2)$.
Let also $\delta_{jk}=\alpha_{jk}-\lambda_{jk}$.

So that with only one matching coordinate between $i_1$ and $i_2$,
$$\eqalign{
Q(i_1,i_2)\indic{|i_1\cap i_2|=1}
&= \Enfa\sum_{k_1,k_2}\delta_{jk_1}\delta_{jk_2} 
\bigl(
\Phi_{jk_1}(X)\Phi_{jk_2}(X)
+
\Phi_{jk_1}(R)\Phi_{jk_2}(R)\bigr)\cr
&-2 \sum_k\delta_{jk}\alpha_{jk}\sum_k \delta_{jk}\lambda_{jk}
 }$$

Again by lemma \ref{dbconcentration} and lemma
\ref{raw moments varphi},
for $X$ or $R$
$$\eqalign{
\Enfa\sum_{k_1,k_2}\delta_{jk_1}\delta_{jk_2}| 
\Phi_{jk_1}(X)\Phi_{jk_2}(X)|
&\leq (4N-3)^d \sum_{k}\delta_{jk}^2 \Enfa\Phi_{jk}(X)^2
\leq C \sum_k\delta_{jk}^2\leq C
 }$$
and since all other terms are bounded by  a constant
not depending on $j$,  by lemma \ref{matching indices}
$(A_n^2)^{-2}\sum_{i_1, i_2} Q(i_1,i_2)\indic{\,|i_1\cap i_2|=1}\leq Cn^{-1}$.

Likewise, the maximum order of $Q(i_1,i_2)\indic{|i_1\cap i_2|=2}$
is $\sum_k [\Enfa \Phi_{jk}(X)^2]^2$, and the corresponding bound is
$2^{jd}n^{-2}$.

\endproof

\proposition Full sample $\hat C_j^2$ risk upper bound\par
\definexref{onesampleC2j}{\the\notitre.\the\noproposition}{proposition}
{\it 
Let $\tilde X=X_1,\ldots,X_n$ be an independent, identically distributed sample of $f_A$.
  Assume that $f$ is compactly supported
 and that $\varphi$ is a Daubechies  $D2N$.
 Let $\Enfa$ be the the expectation relative to the joint distribution
 of the sample $\tilde X$.
Let $\hat C^2_j$ be the plug-in estimator defined in \eqref{defC2j},
Then on $\{2^{jd}<n^2\}$
 $$ \Enfa\Bigl(\hat C^2_j(\tilde X) -C^2_j\Bigr)^2 \leq
 C2^{jd}n^{-1}
  $$ 
}
\proof

$$\eqalign{
\Enfa \bigl[\hat C^2_j -C^2_j\bigr]^2
&\leq \Enfa 3\bigl(\sum_k \hat\alpha_{jk}^2-\alpha_{jk}^2\bigr)^2
+ 3  \bigl(\sum_k \hat\lambda_{jk}^2-\lambda_{jk}^2\bigr)^2
+ 3 \bigl(4\sum_k \hat\alpha_{jk}\hat\lambda_{jk}-
\alpha_{jk}\lambda_{jk}\bigr)^2\cr
}$$

By proposition \ref{2ndmomenthatalpha}
the first term is of the order of $2^{jd}n^{-1}$.
By proposition \ref{2ndmoment crossed one sample}
the two other terms are of the order of  $Cn^{-1}+2^{j}n^{-1}\indic{2^{jd}<n^2}$.

\endproof

As is now shown, the rate of $\hat D^2_j(\tilde X)$ computed on the full sample is slower than the one for $\hat\Delta_j^2(\tilde R, \tilde S)$ 
in the two samples setting. 

The reason  is that we cannot
always apply lemma \ref{dbconcentration} allowing to reduce double sums in $k_1, k_2$ to a sum on the diagonal $k_1=k_2$ for translates of the same $\varphi$ functions. Indeed, when a match between multi indices $i_1$ and $i_2$ involves  terms corresponding to margins, it is not guaranteed that a match on observation numbers also corresponds to a match on margin numbers; that is to say, in the product  $\varphi(X^{\ell_1}-k_1)
\varphi(X^{\ell_2}-k_2)$, only once in a while $\ell_1=\ell_2$; so  most of the time we can say nothing about the support of the product,
and the sum spans many more terms, hence the additional factor $2^j$  in the risk bound
for $\hat D^2_j$ on the full sample.  

\filbreak
\proposition  Risk upper bound, full sample --- $f_A$\par
\definexref{onesampleD2j}{\the\notitre.\the\noproposition}{proposition}

{\it
Let $X_1,\ldots,X_n$ be an independent, identically distributed sample of $f_A$.
  Assume that $f$ is compactly supported
 and that $\varphi$ is a Daubechies  $D2N$.
Let $\hat D^2_j$ be the U-statistic estimator defined in \eqref{defD2j},
Then 
 $$ \Enfa\biggl(\hat D^2_j(\tilde X) -\sum_{k\in\smallZ^d} \delta_{jk}^2\biggr)^2 \leq C2^{jd} n^{-2} + C^\star 2^jn^{-1}$$ 
 with $\delta_{jk}$ the coordinate of $f_A-f_A^\star$ and $C^\star=0$ at independence, when $f_A=f_A^\star$.
 
 }
\proof

$\Enfa \left[\hat D^2_j(\tilde X)-\sum_{k\in\smallZ^d} \delta_{jk}^2\right]^2
=\Enfa [\hat D^2_j(\tilde X)]^2 - \left(\sum_{k\in\smallZ^d} \delta_{jk}^2\right)^2$.

To make $\hat D^2_j(\tilde X)$ look more like the usual U-estimator of
$\int (f-g)^2$ for unrelated $f$ and $g$,
we define for $i\in I_n^{2d+2}$,  the dummy slice variables 
 $Y_i=X_{i^1}$, $V_i=(X_{i^2},\ldots X_{i^{d+1}})$,
$Z_i=X_{i^{d+2}}$, $T_i=(X_{i^{d+3}},\ldots X_{i^{2d+2}})$;
so that $Y_i$ and $Z_i$ have distribution
$P_{f_A}$, $V_i$ and $T_i$ have distribution $P_{f_A^\star}
=P_{\fix(f_A,1)}\ldots P_{\fix(f_A,d)}$ (once canonically projected),
and $Y_i$, $V_i$, $Z_i$, $T_i$ are independent variables
under $\Pnfa$.
Next, for $k\in \Z^d$, define the function $\Lambda_{jk}$ as
$$\eqalign{
\Lambda_{jk}(X_{i^1},\ldots,X_{i^d})&= \varphi_{jk^1}(X_{i^1}^1)\ldots
\varphi_{jk^d}(X_{i^d}^d)\quad \forall i \in \Omega_n^d\cr
\Lambda_{jk}(X_i)=
\Phi_{jk}(X_i)
&=
\varphi_{jk^1}(X_{i}^1)\ldots
\varphi_{jk^d}(X_{i}^d)\quad \forall i\in \Omega_n^1=\{1\ldots,n\}
}\eqdef{conventionlambda}$$
with second line taken as a convention.

So that $\hat D_j^2(\tilde X)$ can be written under
the more friendly form
$$
\hat D^2_j(\tilde X)={1\over A_n^{2d+2}}
\sum_{i\in I_n^{2d+2}}\sum_{k\in \smallZ^d}  \big[\Lambda_{jk}(Y_i)-\Lambda_{jk}(V_i)\big]
\big[\Lambda_{jk}(Z_i)-\Lambda_{jk}(T_i)\big],
$$
with $I_n^m=\Inm$.

Following the friendly notation, let $h_{ik}=\big[\Lambda_{jk}(Y_i)-\Lambda_{jk}(V_i)\big]
\big[\Lambda_{jk}(Z_i)-\Lambda_{jk}(T_i)\big]$  be the kernel of $\hat D_j^2(\tilde X)$
at fixed $k$. Then,
$$ [\hat D^2_j(\tilde X)]^2=|I_n^{2d+2}|^{-2}\sum_{i_1,i_2\in I_n^{2d+2}\times I_n^{2d+2}}
\quad\sum_{k_1,k_2\in \smallZ^d\times \smallZ^d} h_{i_1k_1}h_{i_2k_2}.$$

Consider the partitioning sets 
$M_c=\{i_1,i_2 \in I_n^{2d+2}\times I_n^{2d+2}\colon |i_1\cap i_2|=c\}$,
$c=0\ldots,2d+2$, that is to say the set of pairs with $c$ coordinates in common. Equivalently, 
$M_c$ can be defined as   the set 
$\{i_1,i_2 \in I_n^{2d+2}\times I_n^{2d+2}\colon |i_1\cup i_2|=4d+4-c\}$.

According to the partitioning, with $h_i=\sum_kh_{ik}$, 
$$\Enfa [\hat D^2_j(\tilde X)]^2=|I_n^{2d+2}|^{-2}\sum_{c=0}^{2d+2} \sum_{(i_1,i_2) \in M_c}\Enfa  h_{i_1}h_{i_2}.
$$

Let $\lambda_{jk}=\prodalpha$ and $\delta_{jk}=\alpha_{jk}-\lambda_{jk}$.
\unorderedlist

\li On $M_0$, with no match,
$$\Enfa h_{i_1}h_{i_2}\indic{M_0}=\sum_{k_1,k_2}
(\alpha_{jk_1}-\lambda_{jk_1})^2
(\alpha_{jk_2}-\lambda_{jk_2})^2= \left(\sum_k \delta_{jk}^2\right)^2$$

By lemma \ref{matching indices}, the ratio $|M_0|/|I_n^{2d+2}|$
is lower than $1+Cn^{-2}$. So that

$|I_n^{2d+2}|^{-2}\sum_{M_0} \Enfa h_{i_1}h_{i_2} - 
 \left(\sum_k \delta_{jk}^2\right)^2
 = |I_n^{2d+2}|^{-2} |M_0|\Enfa h_{i_1}h_{i_2}\indic{M_0}
 \leq Cn^{-2}$.

\li On $M_1$, assuming the 
match involves  $Y_{i_1}$ and $Y_{i_2}$,
$$\eqalign{
\Enfa h_{i_1}h_{i_2}\indic{M_1}
&=
\sum_{k_1,k_2}\delta_{jk_1}\delta_{jk_2}
\Enfa
\big(\Phi_{jk_1}(Y_{i_1})-\Lambda_{jk_1}(V_{i_1})\big) 
\big(\Phi_{jk_2}(Y_{i_2})-\Lambda_{jk_2}(V_{i_2}) \big)\cr
&=
\sum_{k_1,k_2}\delta_{jk_1}\delta_{jk_2}
\big(
\Enfa\Phi_{jk_1}(X)\Phi_{jk_2}(X)
 -\lambda_{jk_1}\alpha_{jk_2}
-\delta_{jk_1}\lambda_{jk_2}
\big)\cr
&=\Enfa\left(\sum_{k}\delta_{jk}\Phi_{jk}(X)\right)^2
-C_j^2\sum_k \lambda_{jk}\delta_{jk}
-\left(\sum_k \lambda_{jk}\delta_{jk}\right)
\left(\sum_k \alpha_{jk}\delta_{jk}\right)
}%
\eqdef{symmetricphiphi}
$$ with $C_j^2=\sum_k\delta_{jk}^2$.

Next by \eqref{form2} in lemma \ref{dbconcentration}  for the first line, the double
sum in $k$ under expectation is bounded by a constant times the sum restricted to the diagonal $k_1=k_2$
because of the limited overlapping of translates $\varphi_{jk}$; using also lemma \ref{raw moments varphi},
$$\eqalign{
\Enfa\Bigl(\sum_{k}\delta_{jk}\Phi_{jk}(X)\Bigr)^2
\leq (4N-3)^d\sum_k \delta_{jk}^2\Enfa\Phi_{jk}(X)^2
\leq (4N-3)^d
\sum_k C\delta_{jk}^2.
}$$

Since all other terms in \eqref{symmetricphiphi} are clearly bounded
by a constant not depending on $j$, 
we conclude by symmetry that
$\Enfa h_{i_1}h_{i_2}\indic{M_1}\leq C$ for any match of cardinality 1 between  
narrow slices ($Y_{i_1}$ $Y_{i_2}$ or $Z_{i_1}$ $Z_{i_2}$ or
$Y_{i_1}$ $Z_{i_2}$ or $Z_{i_1}$ $Y_{i_2}$). Moreover
$C=0$ when $f_A=f_A^\ast$ \ie at independence, because of the omnipresence
of $\delta_{jk}$, the coordinate of $f_A-f^\star_A$.

\li On $M_1$, if the match is 
between  $Y_{i_1}$ and $V_{i_2}$, a calculus  as in \eqref{symmetricphiphi} yields,
$$\eqalign{
\Enfa h_{i_1}h_{i_2}\indic{M_1}
&=-\sum_{k_1,k_2}\delta_{jk_1}\delta_{jk_2}
\Enfa\Phi_{jk_1}(Y_{i_1})\Lambda_{jk_2}(V_{i_2})
+C_j^2\sum_k \alpha_{jk}\delta_{jk}
+\left(\sum_k \lambda_{jk}\delta_{jk}\right)^2;\cr
}$$which can also be found from line 2 of \eqref{symmetricphiphi}
using the swap $\Phi_{jk}(Y_{i_2}) \longleftrightarrow -\Lambda_{jk}(V_{i_2})$
and 
$\alpha_{jk} \longleftrightarrow -\lambda_{jk}$.

Next, for some $\ell\in \{1,\ldots,d\}$,
$$\eqalign{
\sum_{k_1,k_2}\delta_{jk_1}\delta_{jk_2}
\Enfa\Phi_{jk_1}(Y_{i_1})\Lambda_{jk_2}(V_{i_2})
&=\sum_{k_1,k_2}\delta_{jk_1}\delta_{jk_2}\lambda_{jk_2}^\spec{d-1}
\Enfa\Phi_{jk_1}(X)\varphi_{jk_2^\ell}(X^\ell)\cr
}$$
 with special notation $\lambda_{jk}^\spec r=\alpha_{jk^1}^{p_1}\ldots \alpha_{jk^d}^{p_d}$
for some  $p_i$, $0\leq p_i\leq r$, $\sum_{i=1}^d p_i=r$.

In the present case
$\Phi_{jk_1}(X)\varphi_{jk_2^\ell}(X^\ell)
=\Phi_{jk_1}(X)\varphi_{jk_2^\ell}(X^\ell)\indic{|k_1^\ell-k_2^\ell|<2N-1}$ does not give any useful restriction of the double sum because 
the coefficient $\alpha_{jk}$ hidden in $\delta_{jk}$ is not guaranteed to factorize under any split of dimension unless $A=I$;  
and lemma \ref{dbconcentration} is useless. This is a difficulty that did not raise in propositions \ref{d+1subsample}
and \ref{2subsample} because we could use the fact that these kind of
terms were estimated over independent samples.

Instead write $\Enfa |\Phi_{jk_1}(X)\varphi_{jk_2}(X^\ell)|
\leq 2^{j\over2}\|\varphi\|_\infty\Enfa|\Phi_{jk_1}(X)|\leq 
C2^{j\over2}2^{-{jd\over2}}$ using lemma \ref{raw moments varphi}. So that when multiplied by 
$\sum_k \delta_{jk}\sum_k\delta_{jk}\lambda_{jk}^\spec{d-1}$, using Meyer's lemma, the final order is $2^j$.

By symmetry,  
for any match of cardinality 1 between a narrow and a wide slice ($Y$ or $T$ or equivalent pairing),
 $\Enfa |h_{i_1}h_{i_2}|\indic{M_1}\leq C2^{j}$, with $C=0$ at independence.

\li On $M_1$, if the match is between $V_{i_1}$ and $V_{i_2}$, 
by symmetry with \eqref{symmetricphiphi} or using the swap defined above,
$$\eqalign{
\Enfa h_{i_1}h_{i_2}\indicsmall{M_1}
&=
\sum_{k_1,k_2}\delta_{jk}\delta_{jk'}
\Enfa\Lambda_{jk}(V_i)\Lambda_{jk'}(V_{i_2})
-C_j^2\sum_k \alpha_{jk}\delta_{jk}
-\left(\sum_k \lambda_{jk}\delta_{jk}\right)
\left(\sum_k \alpha_{jk}\delta_{jk}\right),
}%
$$ 
 and for some not necessarily matching $\ell_1,\ell_2\in \{1,\ldots,d\}$ (\ie lemma \ref{dbconcentration}  not applicable),  
$$\eqalign{
\sum_{k_1,k_2}\delta_{jk_1}\delta_{jk_2}
\Enfa\Lambda_{jk_1}(V_{i_1})\Lambda_{jk_2}(V_{i_2})
&=
\sum_{k_1,k_2}\delta_{jk_1}\delta_{jk_2}
\lambda_{jk_1}^\spec{d-1}\lambda_{jk_2}^\spec{d-1}
\Enfa\varphi_{jk_1^{\ell_1}}(X^{\ell_1})\varphi_{jk_2^{\ell_2}}(X^{\ell_2})\cr
&\leq \Bigl(\sum_{k}\delta_{jk}
\lambda_{jk}^\spec{d-1}\Bigr)^2=C2^j
}$$
with last line using Meyer's lemma, and having reduced the term under expectation to a constant by Cauchy-Schwarz inequality and lemma 
\ref{raw moments varphi}.

And we conclude again that, for any match of cardinality 1
between two wide slices ($V$ or $T$ or equivalent),
$\Enfa h_{i_1}h_{i_2}\indic{M_1}\leq C2^{j}$,
with $C=0$ at independence.

By lemma \ref{matching indices}, the ratio $|M_1|/|I_n^{2d+2}\times I_n^{2d+2} |\approx  n^{-1}$, so in summary, the bound for $M_1$ has the order $C^\star 2^jn^{-1}$, with $C^\star=0$ at independence.

\li On $M_c$, $c=2\ldots2d+2$. 

Fix the pair of indexes $(i_1,i_2)\in I_n^{2d+2}\times I_n^{2d+2}$,
we need to bound a term having the form
$$\eqalign{
Q(i_1,i_2)&=\Enfa\sum_{k_1,k_2} \Lambda_{jk}(R_{i_1})\Lambda_{jk}(S_{i_1})
\Lambda_{jk_2}(R'_{i_2})\Lambda_{jk_2}(S'_{i_2})\cr
}
$$
where both slices $R_{i_1}\neq S_{i_1}$
 unrelated  with both slices $R'_{i_2}\neq S'_{i_2}$
are chosen
among any of the dummy $Y$, $V$, $Z$, $T$.
\unorderedlist

\li {\it Narrow slices only.} For a match spanning four narrow slices exclusively,  that is to say 
$(Y_{i_1}=Y_{i_2}) \cap (Z_{i_1}=Z_{i_2})$ or
$(Y_{i_1}=Z_{i_2})\cap (Z_{i_1}=Y_{i_2})$, a case possible on $M_2$ only, 
the general term of higher order is written
$\sum_{k_1,k_2} \Enfa\Phi_{jk_1}(X)\Phi_{jk_2}(X)\Enfa\Phi_{jk_1}(X)
\Phi_{jk_2}(X)$. By lemma \ref{dbconcentration} this is 
again lower than 
$(4n-3)^d\sum_{k} \left[\Enfa\Phi_{jk}(X)^2\right]^2$,  that is $C2^{jd}$. By lemma \ref{matching indices}, this case thus contributes
to the general bound up to $C2^{jd}n^{-2}$.

Three narrow slices only is not possible and two narrow slices 
correspond to the case $M_1$ treated above.

\li {\it Wide slices only.}  For a match spanning wide slices on $M_c$, $c=2,\ldots 2d$, 
a general term with higher order
is written 
$\sum_{k_1,k_2} \Enfa\Lambda_{jk_1}(V_{i_1})\Lambda_{jk_1}(T_{i_1})\Lambda_{jk_2}(V_{i_2})
\Lambda_{jk_2}(T_{i_2})$, with $|i_1\cap i_2|=c$, 
 (an equivalent is obtained by swapping one V with a T ). Since the slices are wide, it is not possible to distribute expectation any further right now: 
 if $V_{i_1}$ is always independent of $T_{i_1}$, both terms
 may depend on $V_{i_2}$, say.
Also matching coordinates on $i_1$, $i_2$ do not necessarily
correspond to matching dimensions $X^\ell$ of the observation, and then lemma \ref{dbconcentration} 
is not applicable. Instead write,
$$Q(i_1,i_2)=\sum_{k_1,k_2}
\lambda_{jk_1}^\spec{2d-c}
\lambda_{jk_2}^\spec{2d-c}
\Enfa\Big[\Lambda_{jk_1}^\spec c(V_{i_1},T_{i_1} )
\Lambda_{jk_2}^\spec c(V_{i_2},T_{i_2})
\Big],$$
 with 
$\Lambda_{jk}^\spec c(V_i,T_i)$ a product of $c$ independent terms of the form 
$\varphi_{jk^{\ell}}(X^\ell)$ spanning at least one of the slices $V_i$, $T_i$.

By definition of $i_1$ and $i_2$, the product of $2c$ terms under expectation  can be split into $c$ inde\-pendent products of two terms. 
So, using $\Enfa |\varphi_{jk^\ell}(X)^2|\leq C$  on each bi-term, the order at the end is
 $C\bigl(\sum_k \lambda_{jk}^\spec{2d-c}\bigr)^2$;
and using Meyer's lemma,
the bound is then  of the order of $C2^{jc}$.

Finally, using lemma \ref{matching indices} as above,  the contribution of this kind of term to the general bound 
is $\sum_{c=1}^{2d} 2^{jc}n^{-c}$.

On $\{2^{j}< n\}\supset \{2^{jd}< n^2\}\supset \{2^{jd}< n\}$,
 this quantity is bounded by
$C2^{j}n^{-1}<C2^{jd}n^{-2}$ and  
 on $\{2^{j}> n\}$ it is unbounded. 
\li {\it Narrow and wide slices} 
Reusing the general pattern above, with $c_w\leq 2d$ matching coordinates on wide slices and $c_r\leq 2$ on narrow slices
$$Q(i_1,i_2)=\sum_{k_1,k_2}
\lambda_{jk_1}^\spec{2d-c_w}
\lambda_{jk_2}^\spec{2d-c_w}
\alpha_{jk_1}^{2-c_r}
\alpha_{jk_2}^{2-c_r}
\Enfa\Big[\Lambda_{jk_1}^\spec c(Y_{i_1}, V_{i_1}, Z_{i_1}, T_{i_1} )
\Lambda_{jk_2}^\spec c(Y_{i_2}, V_{i_2},Z_{i_2},T_{i_2})
\Big],$$
 with 
$\Lambda_{jk}^\spec c(Y_i, V_i,Z_i,T_i)$ a product of $c$ independent terms of the form 
$\varphi_{jk^{\ell}}(X)$ or $\Phi_{jk}(X)$ spanning at least one of the slices $V_i$, $T_i$ and one of the slices $Y_i$, $Z_i$.
As above, the bracket is a product of independent bi-terms,
each under expectation bounded by some constant C, by lemma \ref{raw moments varphi}, using Cauchy-schwarz inequality if needed.
So this is bounded by
$$
\eqalign{
Q(i_1,i_2)&\leq C\sum_{k_1,k_2}
\lambda_{jk_1}^\spec{2d-c_w}
\lambda_{jk_2}^\spec{2d-c_w}
\alpha_{jk_1}^{2-c_r}
\alpha_{jk_2}^{2-c_r}
=C\Bigl(
\sum_{k}
\lambda_{jk}^\spec{2d-c_w}
\alpha_{jk}^{2-c_r}
\Bigr)^2;
}$$
using Cauchy-Schwarz inequality and Meyer's lemma this is bounded by
$
2^{{j\over2}(c_w-d)}
2^{{jd\over2}(c_r-1)}
$
and, with lemma \ref{matching indices},  the contribution to the general bound on $\{2^j<n^2\}\supset \{2^{jd}<n^2\}$
is 
$$2^{-jd}\sum_{a=1}^2\sum_{b=1}^{2d} 
2^{jb\over2}n^{-b}2^{jda\over2}n^{-a}\indic{2^j<n^2}\leq Cn^{-1}$$

\endunorderedlist

\endunorderedlist

Finally on $\{2^{jd}<n^2\}$, $
 \Enfa \hat B_{j}^2 -\Bigl(\sum_k\delta_{jk}^2\Bigr)^2
 \leq C^\star 2^jn^{-1}+2^{jd}n^{-2}.
$
\endproof

\subsection Implementation issues

The statistic $\hat C^2_j$ is a plug-in estimator; its evaluation uses
in the first place the complete estimation of the density $f_A$ and margins; which takes a computing time of the order of $O(n(2N-1)^d)$ where $N$ is the order of the Daubechies wavelet, and $n$ the number of observations.

 In the second place, the actual contrast is a simple function of the 
$2^{jd}+d2^{j}$ coefficients that estimate density $f_A$ and its margins;  the additional computing time is then in  $O(2^{jd})$.

 One can see here the main numerical drawback of the wavelet contrast in its total formulation --- to be of exponential  complexity in dimension $d$ of the problem; but this is by definition the cost of a condition that guarantees  mutual independence of the components in full generality: $d$ sets  $B_1,\ldots, B_d$ are mutually independent
if $P(B_1\cap\ldots\cap B_d)=PB_1\ldots PB_d$ for each of the $2^d$ choices of indices in $\{1,\ldots,d\}$. 

Complexity in $jd$ drops down to $O(d^22^{2j})$ if one concentrates on a pairwise independence,  like in kernel ICA and related methods,
and in the minimum marginal entropy type method of Miller and Fisher III (2003). Pairwise independence is in fact equivalent to mutual independence in the no noise ICA  model and with at most one Gaussian component (Comon, 1994). The minimization used by

The pairwise algorithm  used by 
Miller et Fisher (2003) consists in searching for the minimum
 in each of the  $C^2_d$ free plans of $\R^d$, applying Jacobi rotations to select a particular  plan. A search in each plan is equivalent to the case  $d=2$, 
where the problem is to find the minimum in $\theta$ of a function on $\R$, for  $\theta\in [0,\pi/2]$. To do so, the simplest could be to try out all  points from $0$ to $\pi/2$ along
a grid, or to use bisection type methods.

 U-statistic estimators of  $C^2_j$ have complexity at minimum in  $O(n^2(2N-1)^{2d})$, that is to say quadratic in $n$ as the method of Tsybakov and Samarov (2002) which also attains parametric rate of convergence;  on the other hand the complexity  in $jd$ is probably lowered since the  contrast can be computed by accumulation, without it be necessary to keep all projection in memory, but only a window whose width depends upon the length of the Daubechies filter.

\section Appendix 1 -- Propositions \par

\proposition 2nd moment of $\sum_k\hat \alpha_{jk}^2$
 about $\sum_k \alpha_{jk}^2$
\par
\definexref{2ndmomenthatalpha}{\the\notitre.\the\noproposition}{proposition}
{\it 

Let $X_1,\ldots,X_n$ be an independent, identically distributed sample of $f$, a
compactly supported function defined on $\R^d$.
  Assume that $\varphi$ is a Daubechies $D2N$.
Let $\hat\alpha_{jk}={1\over n}\sum_{i=1}^n
\varphi_{jk^1}(X^1_{i})\ldots\varphi_{jk^d}(X^d_{i})$,
$k\in \Z^d$.

 Then
$
\Enfa\Bigl(\sum_k\hat\alpha_{jk}^2 -
\sum_k\alpha_{jk}^2\Bigr)^2=C2^{jd}n^{-1}+C2^{2jd}n^{-2}\indic{2^{jd}>n}$

}
\proof

For the mean, using lemma \ref{raw moments varphi},
$$\eqalign{\Enfa \sum_k\hat \alpha_{jk}^2
&= {1\over n^2}\sum_{i_1=i_2}
\sum_k \Enfa\Phi_{jk}(X_{i_1})\Phi_{jk}(X_{i_2})
+ {1\over n^2} \sum_{i_1\neq i_2} \sum_k \alpha_{jk}^2\cr
&= {1\over n}
\sum_k \Phi_{jk}(X_{i})^2
+ {n-1\over n}  \sum_k \alpha_{jk}^2 = \sum_k \alpha_{jk}^2 + O({2^{jd}\over n}).
}$$

For the second moment,   
let $M_c=\{i_1, i_2, i_3, i_4 \in\{1,\ldots,n\}\colon |\{i_1\}\cup
\ldots\cup \{i_4\}|=c\}$.
$$
\eqalign{\Enfa \bigl(\sum_k\hat \alpha_{jk}^2\bigr)^2
&= {1\over n^4}\sum_{c=1}^4\sum_{i_1,\ldots,i_4}
\Enfa\sum_{k_1,k_2} 
 \Phi_{jk_1}(X_{i_1})\Phi_{jk_1}(X_{i_2})\Phi_{jk_2}(X_{i_3})
\Phi_{jk_2}(X_{i_4})\indic{M_c}\cr
}$$
On ${c=1}$, the kernel is equal to
$\sum_{k_1,k_2} \Phi_{jk_1}(X)^2\Phi_{jk_2}(X)^2\leq 
(4N-3)^d
\sum_{k} \Phi_{jk}(X)^4$ by lemma \ref{dbconcentration}.
And by lemma \ref{raw moments varphi}, $\Enfa \sum_{k} \Phi_{jk}(X)^4\leq \sum_k C2^{jd}=C2^{2jd}$.

On ${c=2}$, the kernel takes three generic forms: (a)
$\sum_{k_1,k_2} \Phi_{jk_1}(X)\Phi_{jk_1}(Y)\Phi_{jk_2}(X)\Phi_{jk_2}(Y)$
or\br
 (b) $\sum_{k_1,k_2} \Phi_{jk_1}(X)^2\Phi_{jk_2}(Y)^2$
or (c)
$\sum_{k_1,k_2} \Phi_{jk_1}(X)\Phi_{jk_1}(Y)\Phi_{jk_2}(Y)^2$.
In cases (a) and (c), using lemma \ref{dbconcentration}, the double sum can be reduced to the diagonal $k_1=k_2$.
So using also lemma \ref{raw moments varphi},
$$\eqalign{
(a)\quad& \Enfa |\sum_{k_1,k_2} \Phi_{jk_1}(X)\Phi_{jk_1}(Y)\Phi_{jk_2}(X)\Phi_{jk_2}(Y)|\leq 
\Enfa (4N-3)^d \sum_{k} \Phi_{jk}(X)^2\Phi_{jk}(Y)^2\leq C2^{jd}\cr
(b)\quad & \Enfa \sum_{k_1,k_2} \Phi_{jk_1}(X)^2\Phi_{jk_2}(Y)^2 \leq C2^{2jd}\cr
(c)\quad &
\Enfa |\sum_{k_1,k_2} \Phi_{jk_1}(X)\Phi_{jk_1}(Y)\Phi_{jk_2}(Y)^2|
\leq \Enfa (4N-3)^d \sum_k 
|\Phi_{jk}(X)\Phi_{jk}(Y)^3|\leq C2^{jd}.
}$$
On $c=3$ the only representative form is 
$$\Enfa\sum_{k_1,k_2}  \Phi_{jk_1}(X)\Phi_{jk_1}(Y)\Phi_{jk_2}(Z)^2
= \sum_k \alpha_{jk}^2\sum_k \Enfa \Phi_{jk}(X)^2
\leq C2^{jd},$$
and on $c=4$ the statistic is unbiased equal to $\bigl(\sum_k\alpha_{jk}^2\bigr)^2$ under expectation.

Next, since $|M_4|=A_n^4$ and, using lemma \ref{generalized matching indices},
$|M_c|=O(n^c)$,
$$
\eqalign{\Enfa \bigl(\sum_k\hat \alpha_{jk}^2\bigr)^2
&\leq  A_n^4n^{-4}\bigl(\sum_k\alpha_{jk}^2\bigr)^2
+ C2^{2jd}n^{-3} + n^{-2}2^{2jd} + n^{-1}2^{jd}
\cr
&\leq \bigl(\sum_k\alpha_{jk}^2\bigr)^2 + Cn^{-2}
+ Cn^{-1}2^{jd}\indic{2^{jd}<n} +
Cn^{-2}2^{2jd} \indic{2^{jd}>n} 
}$$ 
with $A_n^4n^{-4}=1-{6\over n}+Cn^{-2}$.

Finally 
$$
\eqalign{\Enfa \bigl(\sum_k\hat \alpha_{jk}^2
-\sum_k\alpha_{jk}^2\bigr)^2
&= \Enfa \bigl(\sum_k\hat \alpha_{jk}^2\bigr)^2
+\bigl(\sum_k\alpha_{jk}^2\bigr)^2 
-2 \Enfa \sum_k\hat \alpha_{jk}^2\sum_k\alpha_{jk}^2\cr 
&\leq  Cn^{-2}
+ Cn^{-1}2^{jd}+
Cn^{-2}2^{2jd} \indic{2^{jd}>n} 
}$$ 

\endproof

\proposition 2nd moment of $\sum_k\hat \lambda_{jk}^2$
 about
$\sum_k \lambda_{jk}^2$
and of $\sum_k\hat \lambda_{jk}\hat\alpha_{jk}$ 
about $\sum_k\lambda_{jk}\alpha_{jk}$ 
\par
\definexref{2ndmoment crossed one sample}{\the\notitre.\the\noproposition}{proposition}
{\it 

Let $X_1,\ldots,X_n$ be an independent, identically distributed sample of $f$, a
compactly supported function defined on $\R^d$.
  Assume that $\varphi$ is a Daubechies $D2N$.
Let $\hat\lambda_{jk}={1\over n^d}\sum_{i=1}^n
\varphi_{jk^1}(X^1_{i})\ldots\sum_{i=1}^n\varphi_{jk^d}(X^d_{i})$,
$k\in \Z^d$.

 Then  on $\{2^{jd}<n^2\}$
$$\eqalign{
\Enfa\Bigl(\sum_k\hat\lambda_{jk}\hat\alpha_{jk} -
\sum_k\lambda_{jk}\alpha_{jk}\Bigr)^2&\leq O(n^{-2})+C{2^j\over n}\cr
\Enfa\Bigl(\sum_k\hat\lambda_{jk}^2 -
\sum_k\lambda_{jk}^2\Bigr)^2& \leq O(n^{-2})+C{2^j\over n}\cr
}$$

}
\proof

$$\eqalign{
\Enfa \Bigl(\sum_{k} \hat\lambda_{jk}^2-\lambda_{jk}^2\Bigr)^2=
 \Enfa \left[\Bigl(\sum_k \hat\lambda_{jk}^2\Bigr)^2
-2\sum_{k}\lambda_{jk}^2\sum_{k} \hat\lambda_{jk}^2
+ \Bigl(\sum_{k} \lambda_{jk}^2\Bigr)^2\right]
}$$

 For $i\in \Omega_n^{2d}$, let $V_i$ be the slice 
$(X^1_{i^1}, X^1_{i^2},\ldots, X^d_{i^{2d-1}},X^d_{i^{2d}} )$. Let
the coordinate-wise kernel function $\Lambda_{jk}$ be given by
 $\Lambda_{jk}(V_i)=\varphi_{jk^1}(X^1_{i^1})\varphi_{jk^1}(X^1_{i^2})\ldots \varphi_{jk^d}(X^d_{i^{2d-1}})\varphi_{jk^d}(X^d_{i^{2d}})$.

Let $|i|$ be the shortcut
notation for 
$|\{i^1\}\cup\ldots\cup \{i^{2d}\}|$. Let $W_n^{2d}=\{i \in \Omega_n^{2d}\colon |i|<2d\} $, that is to say the set of indices with at least one  repeated coordinate.

Then the mean term is written
$$\eqalign{
\Enfa \sum_k\hat \lambda_{jk}^2&=n^{-2d}\sum_{i\in \Omega_n^{2d}}\sum_k \Lambda_{jk}(V_i)\cr
&=n^{-2d}\sum_{W_n^{2d}}\sum_k \Enfa
\Lambda_{jk}(V_i) + A_n^{2d}n^{-2d}\sum_k \lambda_{jk}^2\cr
&=Q_1 +  A_n^{2d}n^{-2d}\theta
}
$$

Let $M_c=\{i\in \Omega_n^{2d}\colon |i|=c\}$ be the set  indices
with $c$ common coordinates. 
So that $Q_1$ is written
$$
Q_1=n^{-2d}\sum_{c=1}^{2d-1}\indic{M_c}\sum_{M_c}\sum_k
 \Enfa
\Lambda_{jk}(V_i)= \sum_k Q_{1jk}
$$ 
 
By lemma \ref{bound kernel products} with lemma parameters $(d=1, m=2d, r=1)$,
$\Enfa|
\Lambda_{jk}(V_i)|\indic{M_c}\leq C2^{{j\over2}(2d-2c)}$ and 
by lemma \ref{generalized matching indices}, $|M_c|=O(n^c)$.
Hence,
$$
Q_{1jk}\leq \sum_{c=1}^{2d-1}n^{-2d+c}C2^{j(d-c)}
=2^{-jd}\sum_{c=1}^{2d-1} C\left({2^j\over n}\right)^{(2d-c)}
$$ 
which on $\{2^{jd}<n\}$ has maximum order $2^{j(1-d)}n^{-1}$ when $d-c$ is minimum \ie $c=2d-1$. 
Finally $|Q_1|\leq \sum_kC2^{j(1-d)}n^{-1} \leq C2^jn^{-1}$.

\medskip 
Next, the second moment about zero is written
$$\eqalign{
\Enfa \Bigl(\sum_k\hat \lambda_{jk}^2\Bigr)^2&=n^{-4d}
\sum_{i_1,i_2\in (\Omega_n^{2d})^2}\sum_{k_1,k_2} \Lambda_{jk_1}(V_{i_1})
\Lambda_{jk_2}(V_{i_2})\cr
&=n^{-4d}\sum_{W_n^{4d}}\sum_{k_1,k_2}
\Enfa\Lambda_{jk_1}(V_{i_1})
\Lambda_{jk_2}(V_{i_2}) 
 + A_n^{4d}n^{-4d}\big(\sum_k \lambda_{jk}^2\Bigr)^2\cr
 &= Q_2 +  A_n^{4d}n^{-4d}\theta^2
}
$$ 
with  
$W_n^{4d}=\{i_1,i_2 \in (\Omega_n^{2d})^2\colon |i_1\cup i_2|<4d\} $, that is to say the set of indices with at least one  repeated coordinate somewhere.

Let this time $M_c=\{i_1,i_2\in (\Omega_n^{2d})^2\colon |i_1\cup i_2|=c\}$ be the set  indices
with overall $c$ common coordinates in $i_1$ and $i_2$. 
So that $Q_2$ is written
$$
Q_2=n^{-4d}\sum_{c=1}^{4d-1}\indic{M_c}\sum_{M_c}\sum_{k_1,k_2}
 \Enfa
\Lambda_{jk_1}(V_{i_1})\Lambda_{jk_2}(V_{i_2})
=\sum_{k_1,k_2}
 Q_{2j_1k_1j_2k_2}
$$ 

By lemma \ref{nomatching across dimensions}, unless $c=1$, it is always possible to find indices $i_1,i_2$
with no match between the observations falling under $k_1$ and those falling under $k_2$, so that there is no way to reduce the double sum in $k_1, k_2$ to a sum on the diagonal using lemma \ref{dbconcentration}.
Note that if $c=1$, $\Enfa\Lambda_{jk}(V_i)\Lambda_{jk}(V_i)
= \Enfa \Phi_{jk}(X)^4$ has order $C2^{jd}$.

So coping with the double sum, by lemma \ref{bound kernel products} with lemma parameters
$(d=1, m=2d, r=2)$,
$\Enfa |\Lambda_{jk}(V_{i_1})\Lambda_{jk}(V_{i_2})|
\leq C2^{{j\over 2}(4d-2c)}$, and again by lemma \ref{generalized matching indices} $|M_c|=O(n^c)$, so
$\Enfa|Q_{2j_1k_1j_2k_2}|
\leq
\sum_{c=1}^{4d-1}n^{c-4d}C2^{{j\over 2}(4d-2c)} 
$, which on $\{2^{jd}<n\}$ has maximum order $2^{j(1-2d)}n^{-1}$ when $c=4d-1$.
Finally, $\Enfa Q_2 \leq \sum_{k_1,k_2}C
2^{j(1-2d)}n^{-1}\leq C2^{j}n^{-1}$.

Putting all together, and 
since
 $A_n^pn^{-p}= 1-{(d+1)(d+2)\over 2n}+O(n^{-2})$,
$$\eqalign{
\Enfa \Bigl(\sum_k \hat\lambda_{jk}^2-\lambda_{jk}^2\Bigr)^2
&= Q_2 + A_n^{4d}n^{-4d}\theta^2 
-2\theta(Q_1 + A_n^{2d}n^{-2d}\theta) + \theta^2\cr
&= Q_2 -2\theta Q_1+
\theta^2(1+ A_n^{4d}n^{-4d}-2  A_n^{2d}n^{-2d})\leq |Q_2| + 2\theta |Q_1| + O(n^{-2})\cr
&\leq  C2^{j}n^{-1}
}
$$

\bigskip
For the cross product, 

As above, 
 for $i\in \Omega_n^{d+1}$, let $V_i$ be the slice 
$(X_{i^0}, X^1_{i^1},\ldots, X^d_{i^{d}} )$. Let
the coordinate-wise kernel function $\Lambda_{jk}$ be given by
 $\Lambda_{jk}(V_i)=\Psi_{jk}(X_{i^0})\psi_{jk^1}(X^1_{i^1})\ldots \psi_{jk^d}(X^d_{i^{d}})$. Let $\theta=\sum_k\alpha_{jk}\lambda_{jk}$.

 Let $W_n^{d+1}=\{i \in \Omega_n^{d+1}\colon |i|<d+1\} $, that is to say the set of indices with at least one  repeated coordinate.

So that,
$\Enfa\sum_k
\hat\alpha_{jk}\hat\lambda_{jk} = Q_1 + A_n^{d+1}n^{-d-1}\theta$
with 
$Q_1=n^{-d-1}\sum_{W_n^{d+1}}\sum_k \Enfa
\Lambda_{jk}(V_i)$ and likewise
$$\Enfa\Bigl(\sum_k
\hat\alpha_{jk}\hat\lambda_{jk}\Bigr)^2=
 Q_2 +  A_n^{2d+2}n^{-2d-2}\theta^2$$
 with 
 $Q_2=n^{-2d-2}\sum_{W_n^{2d+2}}\sum_{k_1,k_2}
\Enfa\Lambda_{jk_1}(V_{i_1})
\Lambda_{jk_2}(V_{i_2})$. And we obtain in the same way,
$$\eqalign{
\Enfa \Bigl(\sum_k \hat\alpha_{jk}\hat\lambda_{jk}-\alpha_{jk}\lambda_{jk}\Bigr)^2
\leq |Q_2| + 2\theta |Q_1| + O(n^{-2})\cr
}
$$

Let $M_c=\{i\in \Omega_n^{d+1}\colon |i|=c\}$ be the set  indices
with $c$ common coordinates. 
So that $Q_1$ is written
$$
Q_1=n^{-d-1}\sum_{c=1}^{d}\indic{M_c}\sum_{M_c}\sum_k
 \Enfa
\Lambda_{jk}(V_i)= \sum_k Q_{1jk}
$$ 
 
By lemma \ref{bound kernel products} with lemma parameters $(m_d=1, m_1=d, r=1)$,
$$\Enfa|
\Lambda_{jk}(V_i)|\indic{M_c}\leq C2^{{jd\over2}(1-2c_d)}2^{{j\over2}(d-2c_1)}$$ with $c_1+c_d=c$, $0\leq c_1\leq d$, $1\leq c_d\leq 1$ and 
by lemma \ref{generalized matching indices}, $|M_c|=O(n^c)$.
Hence,
$$
Q_{1jk}\leq \sum_{c=1}^{d}n^{-d-1+c}C2^{j(d-dc_d-c_1)}
=2^{j(-1+(1-d)c_d)}\sum_{c=1}^{d} C\left({2^j\over n}\right)^{(d+1-c)}
$$ 
which on $\{2^{jd}<n\}$ has maximum order $C2^{j(1-d)}n^{-1}$ when $d+1-c$ is minimum \ie $c=d$. 
Finally $|Q_1|\leq \sum_k C2^{j(1-d)}n^{-1} \leq C2^jn^{-1}$.

\medskip 
Next, as above 
$
Q_2=\sum_{k_1,k_2}
 Q_{2jk_1jk_2}
$,
and again by lemma \ref{nomatching across dimensions}, unless $c=1$, it is always possible to find indices $i_1,i_2$
with no matching coordinates corresponding also to matching dimension
number, so that there is no way to reduce the double sum in $k_1,k_2$ to a sum on the diagonal using lemma \ref{dbconcentration}.

So coping once more with the double sum, by lemma \ref{bound kernel products} with lemma parameters
$(m_d=1, m_1=d, r=2)$,
$\Enfa |\Lambda_{jk}(V_{i_1})\Lambda_{jk}(V_{i_2})|
\leq C2^{{jd\over 2}(2-2c_d)}2^{{j\over 2}(2d-2c_1)}$,
with $c_1+c_d=c$, $1\leq c_d\leq 2$, $0\leq c_1\leq 2d$,  and again by lemma \ref{generalized matching indices} $|M_c|=O(n^c)$, so
$$\Enfa|Q_{2j_1k_1j_2k_2}|
\leq
\sum_{c=1}^{2d+1}n^{c-2d-2}C2^{j(d-dc_d+d-c_1)} 
=2^{j(-2+(1-d)c_d)}\sum_{c=1}^{2d+1} C\left({2^j\over n}\right)^{(2d+2-c)},
$$ which on $\{2^{jd}<n\}$ has maximum order $C2^{-jd}n^{-1}$ when $c=2d+1$.
Then either $c_d=1$, which means that the two terms $\Phi_{jk_1}(X_{i_1})
\Phi_{jk_2}(X_{i_2})$ match on the observation number, in which
case the sum in $k_1,k_2$ can be reduced; either $c_d=2$.
In the first case the order is 
 $\Enfa Q_2 \leq (4N-3)^d\sum_{k}C
2^{-jd}n^{-1}\leq Cn^{-1}$ and in the second case
 $\Enfa Q_2 \leq \sum_{k_1,k_2}C
2^{1-2jd}n^{-1}\leq C2^jn^{-1}$.

\endproof

\proposition  Variance of $\sum_k\hat B^2_j$\par
\definexref{variance hatb2j}{\the\notitre.\the\noproposition}{proposition}
{\it
Let $\{X_1,\ldots, X_n\}$ be an i.i.d. sample with density $f$.  Assume that $f$ is compactly supported 
 and that $\varphi$ is a Daubechies  $D2N$.

Let $\hat B^2_j=\sum_k {1\over A_n^2} \sum_{i\in I_n^2} \Phi_{jk}(X_{i^1})\Phi_{jk}(X_{i^2})$ be the U-statistic estimator of $\sum_k \alpha_{jk}^2$.

Then on $\{2^{jd}<n^2\}$, 
$$\Enfa \Bigl(\hat B^2_j- \sum_k\alpha_{jk}^2\Bigr)^2 \leq Cn^{-1}
+ 2^{jd}n^{-2}
$$
}
\proof

Write that,
$$
 \eqalign{
 \Enfa\bigl[\hat B^2_j(\tilde  X)\bigr]^2&= n^{-2}(n-1)^{-2} 
 \sum_{i_1,i_2\in I_n^2} \sum_{k_1,k_2}
 \Phi_{jk_1}(X_{i^1_1})
 \Phi_{jk_1}(X_{i^2_1})
 \Phi_{jk_2}(X_{i^1_2})
 \Phi_{jk_2}(X_{i^2_2})\cr
}$$
On $M_4=\{i_1, i_2\in I_n^2\colon |i_1\cup i_2|=4\}$, \ie with no match between the two indices, the kernel 
$h_{i_1}h_{i_2}= \sum_{k_1,k_2}
 \Phi_{jk_1}(X_{i^1_1})
 \Phi_{jk_1}(X_{i^2_1})
 \Phi_{jk_2}(X_{i^1_2})
 \Phi_{jk_2}(X_{i^2_2})$ is unbiased, equal under expectation  to 
$(\sum_k \alpha_{jk}^2)^2$. 

On $M_c, c=2,3$, with at least one match between $i_1$ and $i_2$
lemma \ref{dbconcentration} is applicable to reduce the double sum in $k_1, k_2$ and,
 $$
 \eqalign{
 \Enfa h_{i_1}h_{i_2}\indic{M_2\cup M_3}&= 
 \sum_{i_1,i_2\in I_n^2} \sum_{k_1,k_2}
 \Phi_{jk_1}(X_{i^1_1})
 \Phi_{jk_1}(X_{i^2_1})
 \Phi_{jk_2}(X_{i^1_2})
 \Phi_{jk_2}(X_{i^2_2})\indic{M_2\cup M_3}\cr
 &\leq
 \sum_{M_2,M_3} (4N-3)^d\sum_{k}|
 \Phi_{jk}(X_{i^1_1})
 \Phi_{jk}(X_{i^2_1})
 \Phi_{jk}(X_{i^1_2})
 \Phi_{jk}(X_{i^2_2})|\cr
 &\leq
 \sum_{M_2,M_3} C\sum_{k}
 2^{jd(2-|i_1\cup i_2|)}= C
 \sum_{M_2,M_3} 
 2^{{jd}(3-|i_1\cup i_2|)},
 }$$
using lemma \ref{bound kernel products} with parameter $m=2$
and $r=2$ for line 3.

Next,  by lemma \ref{generalized matching indices},
$|M_c|=O(n^c)$ and 
 $|M_4|$ divided by  $(A_n^2)^2$  is more precisely equal to $1-4n^{-1}+Cn^{-2}$. So that
  $$
 \eqalign{
 \Enfa\bigl[\hat B^2_j(\tilde X)\bigr]^2&\leq (1+Cn^{-2})\bigl(\sum_k \alpha_{jk}^2\bigr)^2 +C 
  \sum_{c=2}^3 n^cn^{-4}
 2^{{jd}(3-c)}
 =\bigl(\sum_k \alpha_{jk}^2\bigr)^2 + Cn^{-1}+ C2^{jd}n^{-2}.
 }$$
\endproof

\proposition  Variance of multisample $\prod \sum_k\hat B^2_j(\tilde R^\ell)$\par
\definexref{variance lambda deux}{\the\notitre.\the\noproposition}{proposition}
{\it
Let $\{R^\ell_1,\ldots, R^\ell_n\}$ be an i.i.d. sample of $\fix(f,\ell)$, $\ell=1,\ldots,d$.  Assume that $f$ is compactly supported
 and that $\varphi$ is a Daubechies  $D2N$.
Assume that the $d$ samples are independent.

Let $\hat B^2_j(R^\ell)=\sum_k {1\over A_n^2} \sum_{i\in I_n^2} \Phi_{jk}(R^\ell_{i^1})\Phi_{jk}(R^\ell_{i^2})$ be the U-statistic estimator of $\sum_k \alpha_{jk^\ell}^2$, $\ell=1\ldots d$.

Then on $\{2^{jd}<n^2\}$, 
$$\Enfa \Bigl(\prod_{\ell=1}^d\hat B^2_j(R^\ell)
- \sum_{k^1,\ldots,k^d}\alpha_{jk^1}^2\ldots\alpha_{jk^d}^2\Bigr)^2 \leq Cn^{-1} + C {2^j\over n^2}.
$$
}
\proof

Successive application of $ab-cd=(a-c)b+ (b-d)c$ leads to
$$
a_1\ldots a_d - b_1\ldots b_d=
\sum_{\ell=1}^d (a_\ell-b_\ell)b_1\ldots b_{\ell-1}a_{\ell+1}\ldots a_d.
\eqdef{telescopic}$$ 
So applying \eqref{telescopic},
$$\eqalign{
\sum_k \hat\lambda_{jk}^2-\lambda_{jk}^2
&=\sum_{k^1\ldots k^d}
\hat\alpha_{jk^1}^2\ldots\hat\alpha_{jk^d}^2
-
\alpha_{jk^1}^2\ldots\alpha_{jk^d}^2\cr
&= 
\sum_{k^1\ldots k^d}
\sum_{\ell=1}^d (\hat\alpha_{jk^\ell}^2-\alpha_{jk^\ell}^2)
\alpha_{jk^1}^2\ldots \alpha_{jk^{\ell-1}}^2
 \hat\alpha_{jk^{\ell+1}}^2\ldots\hat\alpha_{jk^d}^2\cr
&= \sum_{\ell=1}^d C
\sum_{k^\ell} (\hat\alpha_{jk^\ell}^2-\alpha_{jk^\ell}^2)
\sum_{k^{\ell+1}} \hat\alpha_{jk^\ell}^2\ldots
\sum_{k^d}\hat\alpha_{jk^d}^2
 }$$
And
$$\eqalign{
\Bigl(\sum_k \hat\lambda_{jk}^2-\lambda_{jk}^2\Bigr)^2
&\leq d \sum_{\ell=1}^d C
\Bigl(\sum_{k^\ell} (\hat\alpha_{jk^\ell}^2-\alpha_{jk^\ell}^2)
\sum_{k^{\ell+1}} \hat\alpha_{jk^\ell}^2\ldots
\sum_{k^d}\hat\alpha_{jk^d}^2\Bigr)^2
 }$$

Label $Q=\Enfa\Bigl(\sum_k \hat\lambda_{jk}^2-\lambda_{jk}^2\Bigr)^2$.

If the $d$ samples are independent, if $2^{jd}<n^2$,  and by proposition 
\ref{variance hatb2j} with parameter $d=1$, 
 $$\eqalign{Q
&\leq \sum_{\ell=1}^{d-1} \left[C
 (Cn^{-1}+{2^{j}\over n^2})\prod_{l=\ell+1}^{d-1}
 (C+Cn^{-1}+{2^{j}\over n^{2}})\right]
 + C(Cn^{-1}+{2^{j}\over n^2})\cr
&\leq C n^{-1}+C{2^{j}\over n^2}
 }$$

\endproof

\proposition  Variance of multi sample $\sum_k\hat\alpha_{jk}\hat\lambda_{jk}$\par
\definexref{indepmultisample}{\the\notitre.\the\noproposition}{proposition}
{\it
Let $\{X_1,\ldots, X_n\}$ be an independent, identically distributed sample of $f_A$.
Let $\{R^\ell_1,\ldots, R^\ell_n\}$ be an independent, identically distributed sample of $\fix(f,\ell)$, $\ell=1,\ldots,d$.  Assume that $f$ is compactly supported
 and that $\varphi$ is a Daubechies  $D2N$.
Assume that the $d+1$ samples are independent. Let
$\Enfa$ be the expectation relative to the joint samples.

Then $$\Enfa \Bigl(\sum_k\hat\alpha_{jk}(\tilde X)\hat\lambda_{jk}(\tilde R^1,\ldots \tilde R^d)
- \sum_k\alpha_{jk}\lambda_{jk}\Bigr)^2 \leq Cn^{-1}\indic{2^{j}<n}
+ C2^{jd}n^{-d-1}\indic{2^j>n}$$
}
\proof

Let $Q=\Enfa\Bigl(\sum_{k\in \smallZ^d} \hat\alpha_{jk}\hat\lambda_{jk}\Bigr)^2$; expanding the statistic,
$$\eqalign{
Q
&= \Enfa\sum_{k_1,k_2}{1\over n^{2d+2}}
\sum_{i\in \Omega_n^{2d+2}} \Phi_{jk_1}(X_{i^1})\Phi_{jk_2}(X_{i^2})
\varphi_{jk_1^1}(R^1_{i^3})\varphi_{jk_2^1}(R^1_{i^4})
\ldots\varphi_{jk_1^d}(R^d_{i^{2d+1}})\varphi_{jk_2^d}(R^d_{i^{2d+2}}).\cr
}$$
By independence of the samples, we only need to consider local constraints
on the coordinates of $i\in\Omega_n^{2d+2}$.

Let $a$ be a subset of $\{0,1,\ldots d\}$.
Let $J_a=\{i\in \Omega_n^{2d+2}\colon
\ell \in a \Rightarrow i^{2\ell+1}= i^{2\ell+2};\;
 \ell \notin a\Rightarrow  i^{2\ell+1}\neq i^{2\ell+2}\}$. It is clear that 
$|J_{a}|=\bigl(n(n-1)\bigr)^{d+1-|a|}n^{|a|}$
and  that
the $J_{a}\,$s define a partition of $\Omega_n^{2d+2}$ when 
$a$ describes  the $2^{d+1}$ subsets of 
$\{0,1,\ldots d\}$. One can check that there  are $C_{d+1}^{c}$ distinct sets $a$ 
such that $|a|=c$, and that $\sum_{c=0}^{d+1} C_{d+1}^c n^c(n(n-1))^{d+1-c}
=n^{d+1}\sum_{c=0}^{d+1} C_{d+1}^c (n-1)^{d+1-c}
=n^{2d+2}$.

On $J_\emptyset$ the kernel is unbiased.
On $J_a$, $0\in a$, with the first two coordinates matching,
the  sum in 
$k_1$, $k_2$ can be reduced to a sum on the diagonal by lemma \ref{dbconcentration}. If $0\notin a$, but some $\ell\in a$ the sum can be reduced only on dimension $\ell$, $k_1^\ell=k_2^\ell$, but to no purpose as will be seen below.

So $Q$ is written $Q=n^{-2d-2}\sum_{a\in {\cal P}(\{0,\ldots,d\})} Q_{0a} 
 +Q_{1a}$, with$$\eqalign{
Q_{0a}
&\leq C_1\sum_{i\in J_a, \, 0\in a}\;\sum_{k\in \smallZ^d}
 \Enfa\Phi_{jk}(X)^2 \Enfa\varphi_{jk^{\ell_1}}(R^{\ell_1})^2\ldots\Enfa\varphi_{jk^{\ell_{|a|-1}}}(R^{\ell_{|a|-1}})^2
 \alpha_{jk^{l_1}}^2\ldots\alpha_{jk^{l_{d-|a|+1}}}^2 
}$$
and
$$\eqalign{
Q_{1a}
&= \sum_{i\in J_a, \, 0\notin a}\;\sum_{k_1,k_2}
 \alpha_{jk_1}\alpha_{jk_2} \Enfa\varphi_{jk_1^{\ell_1}}(R^{\ell_1})
 \varphi_{jk_2^{\ell_1}}(R^{\ell_1})
 \ldots
 \Enfa\varphi_{jk_1^{\ell_{|a|-1}}}(R^{\ell_{|a|-1}})
 \varphi_{jk_2^{\ell_{|a|-1}}}(R^{\ell_{|a|-1}})\cr
 &
 \phantom{
 \sum_{i\in J_a, \, 0\notin a}\;\sum_{k_1,k_2}
 \alpha_{jk_1}\alpha_{jk_2}}
 \alpha_{jk_1^{l_1}}\alpha_{jk_2^{l_1}}\ldots
 \alpha_{jk_1^{l_{d-|a|+1}}}\alpha_{jk_2^{l_{d-|a|+1}}} 
}$$
for some all distinct  $\ell_1,\ldots \ell_{|a|-1}$ and 
$l_1,\ldots l_{d-|a|+1}$ whose union is $\{1,\ldots d\}$
and with $C_1=(4N-3)^d$.
The bound for $Q_{0a}$ is also written
$$\eqalign{
(4N-3)^d
\sum_{i\in J_a, \, 0\in a}
\sum_{k\in \smallZ^d} C
 \lambda_{jk}^\spec{2d-2|a|+2}
}$$
with special notation 
 $\lambda_{jk}^\spec r=\alpha_{jk^1}^{p_1}\ldots \alpha_{jk^d}^{p_d}$
for some integers $p_1,\ldots, p_d$, $0\leq p_i\leq r$ 
with $\sum_{i=1}^d p_i=r$.
And so,  by Meyer's lemma this is also bounded by
$\sum_{i\in J_a, \, 0\in a} C2^{j(|a|-1)}
$.

For $Q_{1a}$ with $|a|\geq 1$, the sum in $k_1,k_2$ could be split
in $k_1^{l_1}\ldots k_1^{l_{d-|a|+1}},k_2^{l_1}\ldots k_2^{l_{d-|a|+1}}$ where no concentration on the diagonal is ensured,
and
$k^{\ell_1}\ldots k^{\ell_{|a|-1}}$ where lemma 
\ref{dbconcentration}
is applicable, but precisely the multidimensional coefficient
$\alpha_{jk}=\alpha_{jk^1\ldots k^d}$ is not guaranteed factorisable under any split, unless
$A=I$. So we simply fall back to   
$$\eqalign{
Q_{1a}
&\leq \sum_{i\in J_a, \, 0\notin a}\;\sum_{k_1,k_2}
 \bigl[\alpha_{jk_1}
 \alpha_{jk_2}\bigr]\bigl[
  \alpha_{jk_1^{l_1}}\alpha_{jk_2^{l_1}}\ldots
 \alpha_{jk_1^{l_{d-|a|+1}}}\alpha_{jk_2^{l_{d-|a|+1}}}\bigr] 
 \bigl[C2^{j\over2}\Enfa|\varphi_{jk_1^\ell}(R^\ell)|\bigr]^{|a|-1}.
 }$$
 This is also written, using
Meyer's lemma at the end,
$$\eqalign{
Q_{1a}
&\leq \sum_{i\in J_a, \, 0\notin a}
\Bigl( \sum_{k}
 \alpha_{jk} \lambda_{jk}^\spec{d-|a|+1} \Bigr)^2
\leq \sum_{i\in J_a, \, 0\notin a}C2^{j(|a|-1)}
}$$

Finally, with $\sum_{i\in J_a}1=|J_a|$ given above, the general bound is written,
$$
Q \leq n^{-2d-2}\left[\sum_{a\neq \emptyset} C2^{j(|a|-1)}n^{d+1}(n-1)^{d+1-|a|}
+ 
n^{d+1}(n-1)^{d+1}\Bigl(\sum_k\alpha_{jk}\lambda_{jk}\Bigr)^2\right]
$$and so
$$\eqalign{
Q-\Bigl(\sum_k\alpha_{jk}\lambda_{jk}\Bigr)^2
&\leq 2^{-j}\sum_{c=1}^{d+1}2^{jc}(n-1)^{-c} + Cn^{-2}\cr
&\leq Cn^{-1}\indic{2^j<n}+ 2^{jd}n^{-d-1}\indic{2^j>n}
}$$
\endproof

\section Appendix 2 -- Lemmas
 
\lemma Property set\par
\definexref{overlap}{\the\notitre.\the\nolemma}{lemma}
{\it
Let $A_1,\ldots,A_r$ be $r$  non empty subsets of a finite set $\Omega$.
Let $J$ be a subset of $\{1,\ldots,r\}$.

Define the property set  
$B_J=\left\{x\in \cup A_j\colon x\in \cap_{j\in J} A_j ;\; x\notin \cup_{j\in J^c} A_j\right\}$, that is
to  say the set of elements belonging exclusively to the sets listed through $J$.
Let $b_{J}=|B_{J}|$ and  $b_\kappa=\sum_{|J|=\kappa} b_{J}$. 

Then 
$\sum_{\kappa=0}^r\sum_{|J|=\kappa}  B_J=\Omega$, 
 and
$$
|A_1| \lor \ldots |A_r|\leq 
\sum_{\kappa=1}^r  b_\kappa=|A_1\cup\ldots A_r|
\leq |A_1| + \ldots |A_r|
=\sum_{\kappa=1}^r  \kappa b_\kappa
$$
 with equality for the right part only if $b_\kappa=0$, $\kappa=2\ldots,r$ \ie if all sets are disjoint, and equality for the left part if one set $A_i$ contains all the others.

}

\proof

It follows from the definition that no two different property sets intersect
and that the union of  property sets 
defines  a partition  of $\cup A_i$, hence a partition of $\Omega$ with the addition of
the missing complementary $\Omega-\cup A_i$ denoted by $B_\emptyset$. 
The $B_J$ are also the atoms of the Boolean algebra generated by $\{A_1,\ldots,A_r, \Omega-\cup A_i\}$ with usual set operations.

With $B_\emptyset$, an overlapping of $r$ sets defines a partition of $\Omega$ with cardinality at most $2^r$; there are $C_r^\kappa$ property sets satisfying $|J|=\kappa$, with $\sum_{\kappa=0}^r C_r^\kappa=2^r$.

\endproof
\lemma Many sets matching  indices\par
\definexref{generalized matching indices}{\the\notitre.\the\nolemma}{lemma}
{\it

Let $m\in \N$, $m\geq 1$. Let $\Omega_n^m$ be the set of indices 
$\{(i^1,\ldots, i^m)\colon i^j\in\N,\, 1\leq i^j\leq n\,\}$.
Let $r \in \N$, $r\geq 1$.
Let $I_n^m=\InminOnm$.

For $i=(i^1,\ldots,i^m)\in \Omega_n^m$, let $\tilde \imath=\cup_{j=1}^m \{i^j\}\subset \{1,\ldots,n\}$ be the set of distinct integers in $i$. 

Then, for some constant $C$ depending on $m$, 
 $$
  \#\Big\{  (i_1,\ldots,i_r)\in \bigl(\Omega_n^m\bigr)^r,\colon\,
   |\, \tilde \imath_1\cup\ldots \cup \tilde \imath_r\,|=a\Big\}
    = O(n^a)
   I\left\{\,|\tilde \imath_1|\lor\ldots\lor|\tilde \imath_r|\leq a\leq mr\right\}
  $$
and in corollary
 $ \#\Big\{(i_1,\ldots,i_r)\in \bigl(I_n^m\bigr)^r\,\colon\, |\,i_1\cup\ldots \cup i_r\,|=a\Big\} = O(n^a) 
  I\left\{m\leq a\leq mr\right\}$.

}

\proof

In the setting introduced by lemma	\ref{overlap},
building the compound $(\tilde\imath_1,\ldots,\tilde\imath_r)$ while keeping track of matching indices is achieved  by drawing 
$b^1_{\{1\}}=|\,\tilde \imath_1|$ integers in the $2^0$--partition $b^0_\emptyset=\{1,\ldots,n\}$ thus constituting $\tilde i_1$, then 
$b^2_{\{1,2\}}+b^2_{\{2\}}=|\tilde\imath_2|$ integers
in the $2^1$--partition $\{ b^1_{\{1\}},\, b^1_\emptyset\}$
thus constituting two subindexes from which to build $\tilde\imath_2$, 
then 
$b^3_{\{1,2,3\}}+b^3_{\{2,3\}}+b^3_{\{1,3\}}+b^3_{\{3\}}=|\,\tilde\imath_3|$ integers
in the $2^2$--partition $\{ b^2_{\{1,2\}},\,b^2_{\{1\}},\,b^2_{\{2\}},\, b^2_\emptyset\}$ thus constituting $2^2$ subindexes from which to
build $\tilde\imath_3$,
and so on, up to  
 $b^r_{\{1,\ldots, r\}}+\ldots +b^r_{\{r\}}=|\,\tilde\imath_r|$ integers 
in the cardinality $2^{r-1}$ partition $\{ b^{r-1}_{\{1,\ldots,r-1\}} \ldots ,\,b^{r-1}_\emptyset\}$ thus constituting $2^{r-1}$ subindexes
from which to build $\tilde\imath_r$.
 
The number of ways to draw the subindexes composing the $r$ indexes is then 
$$
A^{b^1_{\{1\}}}_{b^0_\emptyset}\;
A_{b^1_{\{1\}}}^{b^2_{\{1,2\}}}
A_{b^1_\emptyset}^{b^2_{\{2\}}}
\ldots 
A_{b^{r-1}_{\{1,\ldots,r-1\}}}^{b^r_{\{1,\ldots, r\}}}
\ldots
A_{b^{r-1}_\emptyset}^{b^r_{\{r\}}}
\eqdef{Aforoverlap}$$
  with the nesting property
 $b^j_J=b^{j+1}_J+b^{j+1}_{J\cup \{j+1\}}$ (provided $J$ exists at step $j$) 
and $A_n^m={n!\over(n-m)!}$.

At step $j$, the only property set  with cardinality equivalent to $n$, is $B^{j-1}_\emptyset$, while all others have cardinalities lower than $m$; so picking integers inside these light property sets involve cardinalities at most in $m!$ that go in the constants, 
while the pick in $B^{j-1}_\emptyset$ entails a cardinality $
A_{b^{j-1}_\emptyset}^{b^j_{\{r\}}}=A_{n-|\,\tilde\imath_1\cup\ldots\cup \,\tilde\imath_{j-1}|}^{b^j_{\{r\}}}
\approx n^{b^j_{\{r\}}}$.

Note that,  at step $j-1$, $b^{j-1}_\emptyset=n-|\,\tilde\imath_1\cup\ldots\cup \tilde\imath_{j-1}\,|$, because, at step $j$, $b^j_{\{j\}}$
  designates the number of integers in $\tilde\imath_j$ not matching any previous index $\tilde\imath_1,\ldots,\tilde\imath_{j-1}$; so that also 
  $\sum_{j=1}^r b^j_{\{j\}}=|\,\tilde\imath_i\cup\ldots \cup\, \tilde\imath_r|$; and incidentally $\sum_{J\owns j_0} b^j_J=|\,\tilde\imath_{j_0}|$.

The number of integers picked from the big property set at each step is
$$
A_{b^0_\emptyset}^{b^1_{\{1\}}}
A_{b^1_\emptyset}^{b^2_{\{2\}}}
\ldots
A_{b^{r-1}_\emptyset}^{b^r_{\{r\}}} 
$$with $b^j_\emptyset=n-|\,\tilde\imath_1\cup\ldots\cup \,\tilde\imath_{j-1}|$, $b^0_\emptyset=n$
and 
$\sum_{j=1}^r b^j_{\{j\}}=|\,\tilde\imath_i\cup\ldots \cup \,\tilde\imath_r|$.

For large $n$ this is equivalent  to $n^{|\,\tilde\imath_1\cup\ldots\cup \,\tilde\imath_r|}$.

Having drawn the  subindexes,
building the indexes effectively is a matter of iteratively intermixing two sets of $a$ and $b$ elements; an operation  equivalent  to highlighting $b$ cells in a line of $a+b$ cells, which can be done in $C_{a+b}^{b}$ ways, with $C_n^p=A_n^p/p!$. 

Intermixing the subindexes thus involve cardinalities at most in $m!$, that go in the constant C.

Likewise, passing from $\tilde \imath$ to $i$  involve 
cardinalities at most in $C_m^{|\tilde\imath|}$ and no dependence on $n$.

For the corollary, if $i\in I_n^m$ then  $\tilde\imath=i$ and $|\tilde\imath|=m$. If moreover $i^1<\ldots <i^r$,
the number of ways to draw the subindexes is given by replacing occurrences of '$A$' by '$C$' in \eqref{Aforoverlap},
with $C_n^m={n!\over m!(n-m)!}$, 
 which does not change the order in $n$. Also  
there is only one way to intermix subindexes, because of the ordering constraint.

\endproof

\lemma Two sets matching  indices [Corollary and complement]\par
\definexref{matching indices}{\the\notitre.\the\nolemma}{lemma}
{\it
Let $I_n^m$ be the set of indices
$\{(i^1,\ldots, i^m)\colon i^j\in\N,\, 1\leq i^j\leq n,\, 
i^j\neq i^\ell \hbox{ if } i\neq \ell\}$,
 and let 
 ${I'}_n^m$ be the subset of ${I}_n^m$
such that 
$\{i^1 < \ldots <i^m \}$.

Then for $0\leq b\leq m$,
$$\eqalign{
  \#\Big\{(i_1,i_2)\in I_n^m\times I_n^m\,\colon\, |\,i_1\cap i_2\,|=b\Big\}  &= A_n^m A_m^b A_{n-m}^{m-b} C_m^b = O(n^{2m-b})
 \cr
  \#\Big\{(i_1,i_2)\in {I'}_n^m\times {I'}_n^m\,\colon\, |\,i_1\cap i_2\,|=b\Big\}  &= C_n^m C_m^b C_{n-m}^{m-b} = O(n^{2m-b}) 
 }$$

In corollary, with $P$ (resp. $P'$) the mass probability on $(I_n^m)^2$
(resp. $\,({I'}_n^m)^2$),
 $P(|i_1\cap i_2|= b)\approx
  P'(|i_1\cap i_2|= b) = O(n^{-b})$
and $P(|i_1\cap i_2|= 0)=P'(|i_1\cap i_2|= 0)\leq 1 - m^2 n^{-1} + Cn^{-2}$.
}

\proof

For $i_1,i_2\in I_n^m$, the equivalence $|i_1\cap i_2|=b \Longleftrightarrow |i_1\cup i_2|=2m-b$ gives the link
with the general case of lemma \ref{generalized matching indices}.

Reusing the  pattern of lemma \ref{generalized matching indices} in a particular case: 
there are $A_n^m$ ways to constitute $i_1$, there are $A_m^b$ ways to draw $b$ unordered integers from $i_1$ and
$A_{n-m}^{m-b}$ ways to draw $m-b$ unordered integers from $\{1,\ldots,n\}-i_1$.

To constitute $i_2$, intermixing both subindexes of $b$ and $m-b$ integers is equivalent  to highlighting $b$ cells in a line of $m$ cells;
there are $C_m^b$ ways to do so.
On ${I'}_n^m$, by definition,  having drawn the $b$ then $m-b$ ordered distinct integers,
 intermixing is uniquely determined.

Incidentally, one can check that
$\sum_{b=0}^m A_m^b A_{n-m}^{m-b} C_m^b =A_n^m$,  
and that $\sum_{b=0}^m C_m^b C_{n-m}^{m-b}  =C_n^m$.  

Dividing by $\bigl(A_n^m\bigr)^2$ or $\bigl(C_n^m\bigr)^2$, both equivalent to $n^{2m}$,
gives the probabilities. 
Finally for the special case $b=0$, use the fact that
$$
{A_n^m\over A_{n-c}^m}=(1-{c\over n})\ldots(1-{c\over n-m+1})
\leq (1-{c\over n})^m
$$

\endproof


\lemma Product of $r$ kernels of degree $m$\par
\definexref{bound kernel products}{\the\notitre.\the\nolemma}{lemma}
{\it

Let $r\in \N^\ast$. Let $m\geq 1$. Let $(\echn(X))$ be an independent, identically distributed sample of a random variable on $\R^d$. 
Let  $\Omega_n^m$ be the set of indices
$\{(i^1,\ldots, i^m)\colon i^j\in\N,\, 1\leq i^j\leq n\}$
.

For $i\in \Omega_n^m$,
 define 
$$\eqalign{
a_{ik}&= \Phi_{jk}(X_{i^1})\ldots
\Phi_{jk}(X_{i^m})\cr 
%
%
b_{ik}&=\varphi_{jk}(X_{i^{1}}^{\ell_1})\ldots
\varphi_{jk}(X_{i^{m_1}}^{\ell_{m_1}})
\Phi_{jk}(X_{i^{m_1+1}})\ldots
\Phi_{jk}(X_{i^{m_1+m_d}}).
}$$

Let  $\tilde \imath $ be the set of distinct coordinates in $i$
and let $c=c(\tilde \imath_1,\ldots \tilde\imath_r)
= |\tilde \imath_1\cup\ldots \cup \tilde \imath_r|$ be the overall
number of distinct coordinates in $r$ indices 
$(i_1,\ldots i_r)\in ({\Omega_n^m})^r$.

Then
$$\eqalign{
E^n_f  |a_{i_1k_1}\ldots a_{i_rk_r}|
&\leq C2^{{jd\over2}\left(mr -2c \right)}
\cr
%
E^n_f |b_{i_1k_1}\ldots b_{i_rk_r}|&\leq
C2^{{jd\over2}\left(m_dr -2c_d \right)}
\,2^{{j\over2}\left(m_1r -2c +2c_d\right)}
}$$
with $c_d=c_d(\tilde \imath_1,\ldots \tilde\imath_r)\leq c$ the
fraction of $c$  corresponding to products with at least one $\Phi(X)$ term and $1\leq c_d\leq m_d r$, $0\leq c-c_d\leq m_1 r$,
$1\leq c\leq (m_1+m_2)r$.

}
\proof

Using lemma \ref{overlap}, one can see that the product 
 $a_{i_1k_1}\ldots a_{i_rk_r}$, made of $mr$ terms, can always be split into
 $|\tilde \imath_1\cup \ldots \cup \tilde \imath_r|$ independent products of 
 $c(l)$ dependent terms, $1\leq l\leq |\tilde \imath_1\cup \ldots \cup \tilde \imath_r|$, with $c(l)$ in the range  from
 $|\tilde \imath_1|\lor\ldots\lor|\tilde \imath_r|$ to $mr$
 and $\sum_lc(l)=mr$.

Using lemma \ref{raw moments varphi}, a product of $c(l)$ dependent terms, is bounded under expectation by $C2^{{jd\over 2}(c(l)-2)}$. Accumulating all independent products,  the overall order is 
$C2^{{jd\over 2}(mr-2|\tilde \imath_1\cup \ldots \tilde \imath_r|)}$.

For $b_{i_1k_1}\ldots b_{i_rk_r}$  make the distinction
between groups containing at least one $\Phi(X)$ term and the others containing only $\varphi(X^\ell)$ terms. This splits
the number $|\tilde \imath_1 \cup\ldots\cup \tilde \imath_d|$
into $g_{\Phi,\varphi} + g_\varphi$. Let $c_\varphi(l)$ be the number of $\varphi$ terms in a product of $c(l)$ terms, mixed or not.

On the $g_{\Phi,\varphi}$ groups containing $\Phi$ terms, first bound the product of $c_\varphi(l)$ terms by $C2^{{j\over2}c_\varphi(l)}$, and the remaining terms by $C2^{{jd\over2}(c(l)-c_\varphi(l)-2)}$.
On the $g_{\varphi}$ groups with only $\varphi$ terms, bound the product by $C2^{{j\over2}(c_\varphi(l)-2)}$.

The overall order is then
$$
C2^{{jd\over2}\bigl[\bigl(\sum_{l=1}^{g_{\Phi,\varphi}} c(l)-c_\varphi(l)\bigr) -2g_{\Phi,\varphi}\bigr]}
\;2^{{j\over2}\sum_{l=1}^{g_{\Phi,\varphi}}c_\varphi(l)}
\;2^{{j\over2}\bigl[\bigl(\sum_{l=1}^{g_\varphi}c_\varphi(l)\bigr)-2g_{\varphi}\bigr]}.
$$

The final bound is found using $\sum_{l=1}^{g_\varphi}c_\varphi(l)+
\sum_{l=1}^{g_{\Phi,\varphi}}c_\varphi(l)=m_1r$
and 
$\sum_{l=1}^{g_{\Phi,\varphi}} c(l)-c_\varphi(l)=m_dr$.

Rename $c_d=g_{\Phi,\varphi}$ and $c-c_d=g_{\varphi}$.

As for the constraints, in the product of $(m_1+m_d)r$ terms,
it is clear that
$\Phi$ terms have to be found somewhere, so $c_d\geq 1$, which also 
implies that $c-c_d=0$ when $c=1$ (in this case there are no independent group with only $\phi$ terms, but only one big group with all indices equal).
Otherwise $c_d\leq m_dr$ and $c-c_d\leq m_1r$ since there are no more
that this numbers of $\Phi$ and $\phi$ terms in the overall product.

\endproof

\lemma Meyer\par
\definexref{meyerlemma}{\the\notitre.\the\nolemma}{lemma}
{
\it
Let $V_j,, j\in \Z$ an $r$-regular multiresolution analysis of $L_2(\R^n)$
and let $\varphi\in V_0$ be the father wavelet.

There exist two constant $c_2 > c_1 >0$ such 
that for all $p\in [1,\,+\infty]$
and for all finite sum $f(x)=\sum_k \alpha(k)\varphi_{jk}(x)$ one has,

$$c_1\|f\|_p
\leq 2^{jd({1\over2}-{1\over p})}\left(\sum_k |\alpha(k)|^p\right)^{1\over p}
\leq c_2\|f\|_p
$$

}
\proof 
See Meyer (1997)

We use the bound under a special form.

First note that if $f\in B_{sp\infty}$, 
$\|f\|_{sp\infty} = \|P_jf\|_p + \sup_j 2^{js}\|f-P_jf\|_p$
so that $\|f-P_jf\|_p \leq C\|f\|_{sp\infty}2^{-js}$.
So  using \eqref{meyer},
$$\eqalign{
\sum_k |\alpha_{jk}|^p
&\leq C2^{jd(1-p/2)}\|P_jf\|_p^p\leq C2^{jd(1-p/2)}2^{p-1} \big(\,\|f\|_p^p + \|f-P_jf\|_p^p\big)
\cr
&\leq C2^{jd(1-p/2)}2^{p-1} \big(\,\|f\|_p^p + C\|f\|_{sp\infty}^p 2^{-jps}\big)\cr
&\leq C2^{jd(1-p/2)}\|f\|_{sp\infty}^p.
}$$

When applying the lemma to special coefficient
 $\lambda_{jk}^\spec r=\alpha_{jk^1}^{p_1}\ldots \alpha_{jk^d}^{p_d}$
for some integers $p_1,\ldots, p_d$, $0\leq p_i\leq r$ 
with $\sum_{i=1}^d p_i=r$,
we use 
$$ \eqalign{
\sum_{k\in \smallZ^d} |\lambda_{jk}^\spec r|&=
\sum_{k^1\in \smallZ}|\alpha_{jk^1}^{p_1}|\ldots \sum_{k^d\in \smallZ}|\alpha_{jk^d}^{p_d}|\cr
&\leq C2^{{j\over2}(2-p_1)}\|\fix(f,1)\|_{sp_1\infty}^{p_1}
\ldots2^{{j\over2}(2-p_d)}
\|\fix(f,d)\|_{sp_d\infty}^{p_d}\cr
&\leq C2^{{j\over2}(2d-r)}\|\max_\ell \fix(f,\ell)\|_{sr\infty}^r
}$$
so that even if some $p_\ell$ was zero, the result is a $2^j$, 
 which returns the effect 
 of $\sum_{k^\ell}1$.

\endproof
\filbreak
\lemma Path of non matching dimension numbers\par
\definexref{nomatching across dimensions}{\the\notitre.\the\nolemma}{lemma}
{\it

Let $r\in \N$, $r\geq 2$. Let $\Omega_n^m=\Onm$. For $i\in \Omega_n^d$, let $\Lambda_{jk}(V_i)=\varphi_{jk}(X^1_{i^1})\ldots
\varphi_{jk}(X^d_{i^d})$. Let $\tilde\imath$ be the set of distinct
coordinates of $i$.

In the product
$$
\Bigl(\sum_j\sum_k{1\over n^d}\sum_{i\in\Omega_n^d}\Lambda_{jk}(V_i)\Bigr)^r
= {1\over n^{dr}}
\sum_{i_1,\ldots,i_r\in(\Omega_n^d)^r}
\sum_{j_1\ldots j_r}\sum_{k_1\ldots,k_r}
\Lambda_{j_1k_1}(V_{i_1})\ldots\Lambda_{j_rk_r}(V_{i_r})
$$
unless $|\tilde\imath_1\cup\ldots\cup\tilde\imath_r|<r$,
it is always possible to find indices
 $(i_1,\ldots, i_r)$ such that no two functions $\varphi_{jk}$
 $\varphi_{jk'}$ match on observation number.

}
\proof

Let $c=|\tilde\imath_1\cup\ldots\cup\tilde\imath_r|$.
For $1\leq \ell\leq n$, let $\ell^{\otimes d}=(\ell,\ldots,\ell)
\in \Omega_n^d$.

With $r$ buckets of width $d$ defined by the extent of each
index $k_1\ldots,k_r$, and
only $c<r$ distinct observation numbers, once $c$ buckets
have been stuffed with terms $V_{\ell^{\otimes d}}$,
some already used observation number must be reused
in order to fill in the remaining $r-c$ buckets. So that $r-c$
buckets will match on dimension and observation number allowing
to reduce the sum to only $c$ distinct buckets.

Once $c>r$, starting with a configuration
using $V_{\ell_1^{\otimes d}},\ldots V_{\ell_r^{\otimes d}}$
we can always use additional observation numbers to fragment
further the $\ell^{\otimes d}$ terms, which preserves
the empty intersection between buckets.

\endproof

\lemma Daubechies wavelet concentration property\par
\definexref{dbconcentration}{\the\notitre.\the\nolemma}{lemma}
{
\it
Let $r\in \N$, $r\geq 1$. Let $\varphi$ be the scaling function of a Daubechies wavelet $D2N$.
Let $h_{k}$ be the function on $\R^{m}$
defined as a product of translations of $\varphi$
 $$ h_{k}(x_1,\ldots,x_m)=  
 \varphi(x_{1}-k^1)
 \ldots\varphi(x_{m}-k^m),$$
with $k=(k^1,\ldots, k^m)\in \Z^m$.

Then for a Haar wavelet 
$\big[\sum_k h_{k}(x_1,\ldots,x_m)\big]^r
= \sum_{k}h_{k}(x_1,\ldots,x_m)^r$.

For any D2N, 
$$
\bigg(\sum_k |h_{k}(x_1, \ldots x_m)|\bigg)^r\leq 
(4N-3)^{m(r-1)} \sum_k  |h_{k}(x_1, \ldots x_m)|^{r}
\eqdef{form2}$$

}
\proof

With a Daubechies Wavelet $D2N$, whose support
is $[0,2N-1]$ with $\varphi(0)=\varphi(2N-1)=0$ (except for Haar where 
$\varphi(0)=1$), one has the relation 
$$x\mapsto\varphi(x-k)\varphi(x-\ell)=0,\quad \hbox{ for }  |\ell-k|\geq 2N-1;$$
when  $k$ is fixed,
the cardinal of the set $|\ell -k| < 2N-1$ is equal to  $(4N-3)$.

So that, with $k_1,\ldots k_r$ denoting $r$ independent multi-index,
$$\eqalign{
 \big(\sum_k h_{k}\big)^r
= \sum_{k_1}\sum_{k_2\ldots k_r} h_{k_1}\ldots h_{k_r}
\,I( \Delta)
}
$$ with  $\Delta=\{
 |k^{\ell_1}_{i_1}-k^{\ell_2}_{i_2}| < (2N-1);\; i_1,\,i_2=1\ldots r;\; \ell_1,\ell_2=1\ldots m \}$.
  Once $k_1$ say, is fixed,
  the cardinal of $\Delta$ is not greater than $(4N-3)^{m(r-1)}$ and
  is exactly equal to $1$ for Haar, when all $k_1=\ldots=k_r$.

For any Daubechies wavelet, and $r\geq 1$,
using the inequality
 $(|h_{k_1}|^r\ldots |h_{k_r}|^r)^{1\over r} \leq {1\over r}\sum_i |h_{k_i}|^r$,
$$\eqalign{
 \big(\sum_k |h_{k}|\big)^r&\leq
\sum_{k_1,\ldots, k_r} {1\over r} \big(\,  |h_{k_1}|^{r}+ \ldots + 
|h_{k_r}|^{r}\big)\indic\Delta\cr
&= {1\over r} \biggl[\,\sum_{k_1,\ldots, k_r} |h_{k_1}|^{r}\indic\Delta	+ \ldots + \sum_{k_1,\ldots, k_r}
|h_{k_r}|^{r}\indic\Delta \biggr]\cr
 &\leq 
(4N-3)^{m(r-1)} \sum_k  |h_{k}|^{r},
}$$
\endproof

\lemma  $r$th order moment of $\Phi_{jk}$\par
\definexref{raw moments varphi}{\the\notitre.\the\nolemma}{lemma}
{\it

Let $X$ be random variables on $\R^d$ with density $f$. Let $\Phi$ be the tensorial scaling function of an MRA of $L_2(\R^d)$.
Let $\alpha_{jk}=E_f\Phi_{jk}(X)$. Then for $r\in \N^\ast$,
$$ 
E_f|\Phi_{jk}(X)-\alpha_{jk} |^r
\leq 2^r E_f|\Phi_{jk}(X)|^r\leq 2^r
 2^{jd({r\over2}-1)} \|f\|_\infty \|\Phi\|_r^r.
$$

If $\Phi$ is the Haar tensorial wavelet then also
$
E_f\,\Phi_{jk}(X)^r \leq 2^{jd({r\over2}-{1\over2})} \alpha_{jk}
$.

}
\proof

For the left part of the inequality,
$\Bigl(E_f|\Phi_{jk}(X)-\alpha_{jk}|^r\Bigr)^{1\over r}
\leq  
\Bigl(E_f|\Phi_{jk}(X)|^r\Bigr)^{1\over r}
+E_f|\Phi_{jk}(X)|$, and also
$E_f|\Phi_{jk}(X)|\leq \Bigl(E_f|\Phi_{jk}(X)|^r\Bigr)^{1\over r}
 \Bigl(E_f 1\Bigr)^{r-1\over r}$.

For the right part,
$E_f|\Phi_{jk}(X)|^r = 2^{jdr/2}\int |\Phi(2^jx-k)|^rf(x)dx
\leq 2^{jd({r\over2}-1)} \|f\|_\infty \|\Phi\|_r^r$.

Or also if $\Phi$ is positive, 
$$\eqalign{
E_f \Phi_{jk}(X)^r
 &= 2^{{jd\over 2}(r-1)}\int \Phi(2^jx-k)^{r-1}
 \Phi_{jk}(x)f(x)dx\cr
 &\leq 2^{{jd\over 2}(r-1)} \|\Phi\|_\infty^{r-1} \alpha_{jk}.
}$$

\endproof

\filbreak
\section  References

\parskip=10pt

%
%
%


(Bach \& Jordan, 2002)
M.~I.~Jordan F.~R.~Bach.
 Kernel independent component analysis.
 {\it J. of Machine Learning Research}, 3:1--48, 2002.

(Barbedor, 2005)
P. Barbedor.
 Independent components analysis by wavelets.
 {\it Technical report, LPMA Universit\'e Paris 7}, PMA-995, 2005.

(Bell \& Sejnowski, 1995) A.~J.~Bell. T.J.~Sejnowski 
 A non linear information maximization algorithm that performs blind
  separation.
 {\it Advances in neural information processing systems}, 1995.

(Bergh \& L\"ofstr\"om, 1976)
J.~Bergh and J.~L{\"o}str\"om.
 {\it Interpolation spaces}.
 Springer, Berlin, 1976.

(Bickel \& Ritov, 1988)
P. J. Bickel and Y. Ritov. Estimating integrated squared density derivatives: sharp best order of convergence estimates. {\it Sankya Ser} A50, 381-393

(Butucea \& Tribouley, 2006)
C.~Butucea and K.~Tribouley. Nonparametric homogeneity tests.
{\it Journal of Statistical Planning and Inference}, 136(2006), 597--639.

(Birg\'e \& Massart, 1995) L.~Birg\'e and P.~Massart.
Estimation of integral functionals of a density.
{\it Annals of Statistics} 23(1995), 11-29


(Cardoso, 1999) 
J.F. Cardoso.
 High-order contrasts for independent component analysis.
 {\it Neural computations 11}, pages 157--192, 1999.

(Comon, 1994) 
P.~Comon.
 Independent component analysis, a new concept ?
 {\it Signal processing}, 1994.

(Daubechies, 1992)
Ingrid Daubechies.
 {\it Ten lectures on wavelets}.
 SIAM, 1992.

(Devore \& Lorentz, 1993)
R.~Devore, G.~Lorentz. {\it Constructive approximation}. Springer-Verlag,
1993.

(Donoho et al., 1996)
G.~Kerkyacharian D.L.~Donoho, I.M.~Johnstone and D.~Picard.
 Density estimation by wavelet thresholding.
 {\it Annals of statistics}, 1996.

(Gretton et al. 2003)
Alex~Smola Arthur~Gretton, Ralf~Herbrich.
 The kernel mutual information.
 Technical report, Max Planck Institute for Biological Cybernetics,
  April 2003.

(Gretton et al. 2004)
A.~Gretton,  O.~Bousquet, A.~Smola, B.~Sch\"o{}lkopf.
 Measuring statistical dependence with Hilbert-Schmidt norms.
 Technical report, Max Planck Institute for Biological Cybernetics,
  October 2004.

(H\"ardle et al., 1998)
 Wolfgang H{\"
  a}rdle, G\'erard~Kerkyacharian, Dominique~Picard and 
 Alexander~Tsybakov.
 {\it Wavelets, approximation and statistical applications}.
 Springer, 1998.

(Hyvar\"\i nen et al. 2001)
A.~Hyvar\"\i nen, J.~Karhunen. E.~Oja
{\it Independent component analysis}.
 Inter Wiley Science, 2001.

(Hyvarinen \& Oja, 1997)
A.~Hyvarinen and E.~Oja.
 A fast fixed-point algorithm for independent component analysis.
 {\it Neural computation}, 1997.

(Kerkyacharian \& Picard, 1992)
G{\'e}rard~Kerkyacharian Dominique~Picard.
 Density estimation in {B}esov spaces.
 {\it Statistics and Probability Letters}, 13:15--24, 1992.

(Kerkyacharian \& Picard, 1996)
G{\'e}rard~Kerkyacharian Dominique~Picard.
 Estimating non quadratic functionals of a density using Haar wavelets.
 {\it Annals of Statistics}, 24(1996), 485-507.

(Laurent, 1996) B. Laurent. Efficient estimation of a quadratic
functional of a density. {\it Annals of Statistics} 24 (1996), 659 -681.

(Miller and Fischer III, 2003)
E.G. Miller and J.W. Fischer III. ICA using spacings estimates
of entropy. {\it proceedings of the fourth international
symposium on ICA and BSS}. 2003.

(Meyer, 1997)
Yves Meyer.
 {\it Ondelettes et op\'erateurs}.
 Hermann, 1997.

%

(Nikol'ski$\breve{\hbox{\i}}$, 1975) S.M. 
Nikol'ski$\breve{\hbox{\i}}$. Approximation of functions of several variables and imbedding theorems. {\it Springer Verlag, 1975}.

(Peetre, 1975)
Peetre, J.  
New Thoughts on Besov Spaces. Dept. Mathematics, Duke Univ, 1975.
%

%

(Rosenblatt, 1975)
M.~Rosenblatt.
 A quadratic measure of deviation of two-dimensional density estimates
  and a test for independence.
 {\it Annals of Statistics}, 3:1--14, 1975.

(Rosenthal, 1972)
Rosenthal, H. P. 
On the span in lp of sequences of independent random variables. {\it Israel
J. Math}. 8 273--303, 1972.

(Serfling, 1980)
Robert~J. Serfling.
 {\it Approximation theorems of mathematical statistics}.
 Wiley, 1980.

%


(Tsybakov \& Samarov,  2004) 
A.~Tsybakov A.~Samarov.
Nonparametric independent component analysis.
 {\it Bernouilli}, 10:565--582, 2004.

(Tribouley, 2000) K. Tribouley, Adaptive estimation of integrated
functionals. {\it mathematical methods of statistics} 9(2000) p19-38.

(Triebel, 1992)
Triebel, H. Theory of Function Spaces 2. Birkh\"auser, Basel, 1992

{\bf Programs and other runs available at \tt http://www.proba.jussieu.fr/pageperso/barbedor}

\bye